\documentclass[10pt,a4paper,twoside]{amsart}

\usepackage{amsmath, amssymb, amsthm}
\usepackage{amsaddr}
\usepackage{siunitx}
\usepackage[utf8]{inputenc}
\usepackage[T1]{fontenc}
\usepackage{lmodern}
\usepackage{color}
\usepackage{geometry}
\usepackage{chngcntr}
\usepackage{graphicx}
\usepackage{float}
\usepackage{subfig}
\usepackage[demo,abs]{overpic}
\graphicspath{{./figs/}}
\usepackage[normalem]{ulem}
\usepackage{url}

\newcommand{\TP}{\mathrm{TP}}
\newcommand{\tp}{\mathrm{tp}}

\def\R{\mathbb{R}}

\def\cA{\mathcal{A}}
\def\cB{\mathcal{B}}
\def\cC{\mathcal{C}}
\def\cD{\mathcal{D}}

\def\cF{\mathcal{F}}

\def\cI{\mathcal{I}}

\def\cN{\mathcal{N}}

\def\cQ{\mathcal{Q}}

\def\cS{\mathcal{S}}
\def\cT{\mathcal{T}}

\def\a{\alpha}
\def\b{\beta}

\def\d{\delta}

\def\t{\theta}
\def\T{\Theta}

\def\o{\omega}
\def\veps{\varepsilon}

\def\vphi{\varphi}

\def\GD{{\Gamma_{\rm D}}}

\def\p{\partial}

\def\Th{\mathcal{T}_h}
\def\Nh{\mathcal{N}_h}

\def\DD{{\rm D}}

\def\Id{I}

\def\tx{\widetilde{x}}
\def\ty{\widetilde{y}}

\newcommand{\dv}[1]{\,{\mathrm d}#1}

\DeclareMathOperator{\diver}{div}

\DeclareMathOperator{\diam}{diam}

\DeclareMathOperator{\sign}{sign}

\let\oldmarginpar\marginpar
\renewcommand\marginpar[1]{
  \oldmarginpar[\raggedleft\footnotesize #1]
  {\raggedright\footnotesize #1}}

\theoremstyle{definition}
\newtheorem{definition}{Definition}
\newtheorem{remark}[definition]{Remark}

\newtheorem{example}[definition]{Example}

\theoremstyle{plain}

\newtheorem{algorithm}[definition]{Algorithm}
\numberwithin{definition}{section}
\numberwithin{equation}{section}

\begin{document}
\title{Simulating Self-Avoiding Isometric Plate Bending}
\author{Sören Bartels, Frank Meyer, Christian Palus}
\address{Department of Applied Mathematics, University of Freiburg}
\date{\today}
\begin{abstract}
  Inspired by recent results on self-avoiding inextensible curves, we propose and experimentally investigate a numerical method for simulating isometric plate bending without self-intersections. We consider a nonlinear two-dimensional Kirchhoff plate model which is augmented via addition of a tangent-point energy. The resulting continuous model energy is finite if and only if the corresponding deformation is injective, i.\,e. neither includes self-intersections nor self-contact. We propose a finite element method method based on discrete Kirchhoff triangles for the spatial discretization and employ a semi-implicit gradient descent scheme for the minimization of the discretized energy functional. Practical properties of the proposed method are illustrated with numerous numerical simulations, exploring the model behavior in different settings and demonstrating that our method is capable of preventing non-injective deformations.
\end{abstract}
\maketitle

\section{Introduction}

\subsection{Motivation and outline}
The development and utilization of new materials imply the demand for adequate mathematical models to describe resulting material behavior.
Recently, the vast number of technical applications (e.\,g. \cite{ScEb01, SmInLu95, StPuIo11, SuShMi94, ChSvGeEA16}) of thin structures composed of bilayer polymers sparked interest among engineers and applied mathematicians alike and resulted in the formulation and investigation of mathematical models for describing large bending deformations of plates in general \cite{Friesecke02, Friesecke02b}, and the prestress-induced bending of bilayer plates in particular \cite{Schmidt07, Schmidt07b}. 
Consequently, numerical methods for simulating material behavior based on simple plate models \cite{Bartels13, BonNocNto19}, as well as the bilayer models \cite{BaBoMuNo18, BaBoNo17, BaPa20, BonNocNto20} have been formulated and investigated. However, these methods have in common one practically relevant shortcoming: They are not able to recognize self-contact. Especially in the case of bilayer plates, where large deformations are a key feature, this is a practical restriction as the resulting numerical solutions often exhibit severe self-intersections and, thus, cannot coincide with physical observations.

Coming up with a precise mathematical description of the intuitively simple notions of \emph{self-avoidance} and self-contact is a challenging task, even more so, if the model should be amenable to numerical simulation. Due to its practical relevance, both the investigation of theoretical concepts as well as the implementation of software packages have been the topic of many scientific investigations in the past and present, see e.\,g.~\cite{CiaNec87, GonMadSchMos02, Palmer19, PalHea17, Schumacher20, StrMos06, StrMos07}.
Motivated by the results on isotopy-class-preserving closed curves in \cite{BaRe21, BaReRi18}, in this paper we aim to propose and experimentally investigate a \emph{self-avoiding} plate model which tries to overcome the above mentioned shortcomings. 
The outline is as follows: In the remainder of Section~1 we introduce the two-dimensional bending model as well as the self-repulsive tangent-point potential. 
In Section~2 we gather the preliminaries that are needed for the formulation of our proposed method, which we present in Section~3.
In Section~4 we briefly present three different techniques that may be employed to improve efficiency of implementations.
Section~5 contains several numerical experiments in which we examine different aspects and practical properties of the discretized model and our algorithm for its simulation, thereby providing experimental justification of our method.
We conclude the paper with a short summary of our observations in Section~6.

\subsection{Self-avoiding plate model}
We consider a two-dimensional nonlinear Kirchhoff plate model in the bending regime characterized by cubic scaling of a three-dimensional elastic energy with respect to the plate thickness: we describe a deformed plate by the deformation $y \colon \o \to \R^3$ of its mid-plane with flat reference configuration $\o \subset \R^2$ and only consider isometric deformations which comply with the \emph{isometry constraint} $[\nabla y]^\top \nabla y = \Id_2$. For a given appropriately scaled body force $f \colon \o \to \R^3$ and plates with thickness $h \ll \diam{\o}$, the corresponding elastic bending energy
\[
E_\mathrm{Ki}[y] = \frac 1 2 \int_\o |D^2 y|^2 \dv{x} -\int_\o f\cdot y \dv{x}
\]
has been rigorously derived from three-dimensional elasticity in \cite{Friesecke02, Friesecke02b}. In the case of bilayer plates, which are manufactured from compound materials consisting of two layers with slightly different material properties, the limiting 2D energy functional has been rigorously justified in \cite{Schmidt07, Schmidt07b} and is given by
\[
E_\mathrm{bil}[y] = \frac{1}{2}\int_\o |\mathit{II}(y) - \a \Id_2|^2 \dv{x} - \int_\o f\cdot y \dv{x},
\]
where $\mathit{II(y)}$ denotes the second fundamental form of the parametrized surface defined by the deformation $y$ and $\alpha > 0$ is a parameter corresponding to a homogeneous material mismatch between the layers. Sensible boundary conditions for the minimization of these energies model a clamping of the plate on part of its boundary, i.\,e. $y=y_\DD$ and $\nabla y=\phi_\DD$ on a subset $\GD \subset \p \o$ with positive length.

\begin{figure}[t]\centering
  \includegraphics[height=4.4cm]{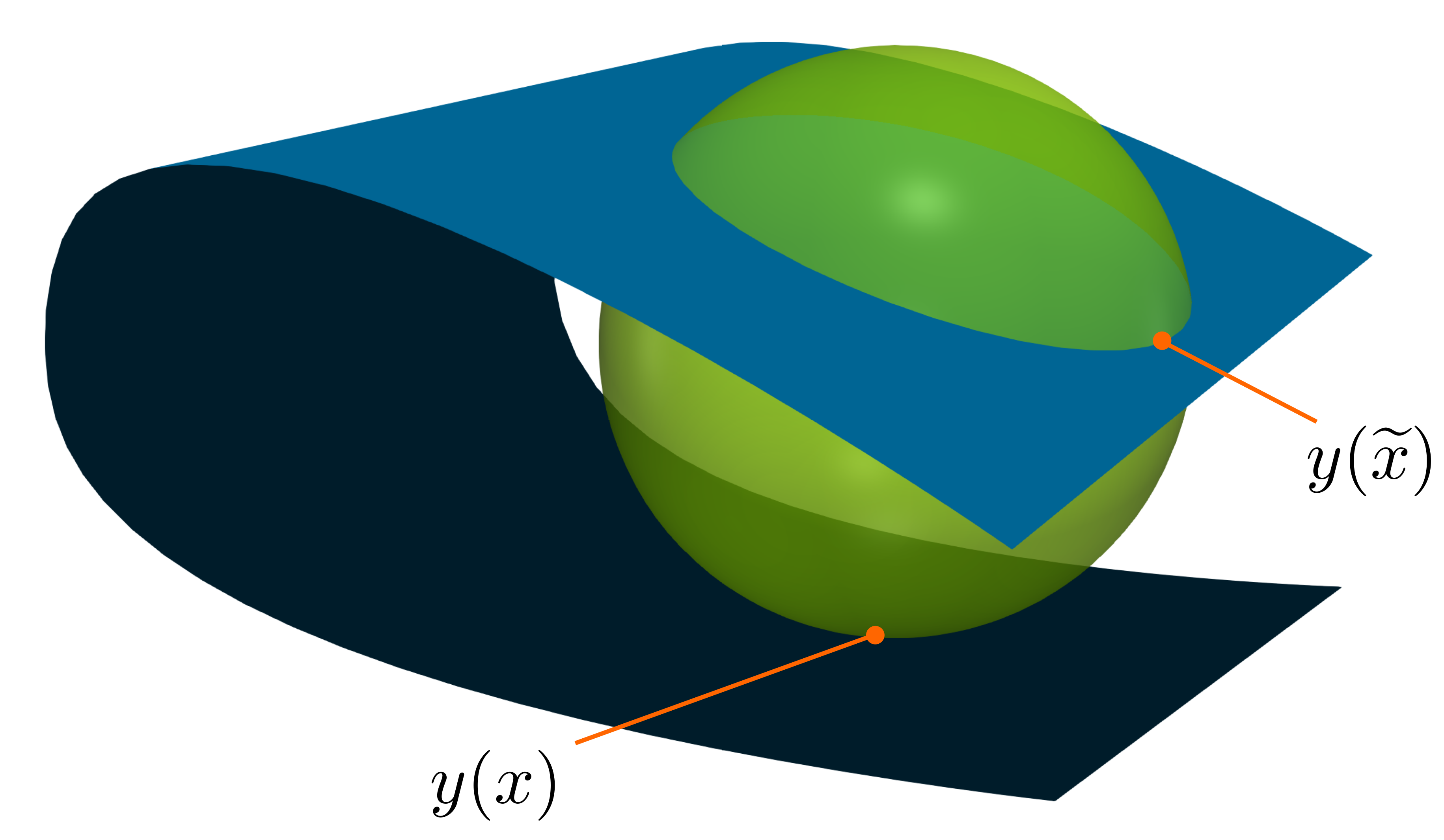}
  \caption{The tangent-point radius $r = r[y](x,\tx)$ is defined by a sphere which is tangent to the surface in $y(x)$ and which intersects the surface in $y(\widetilde{x})$.\label{fig:radius}}
\end{figure}

The existence of minimizers for the resulting constrained minimization problems can be established by means of the direct method in the calculus of variations. In general, however, the minimizers cannot be expected to be one-to-one and thus may exhibit non-physical properties in the form of self-intersections. This is particularly evident in the case of bilayer plates, but can also easily be observed in the single-layer case when compressive boundary conditions are imposed. As a remedy, we consider for $q \ge2$ the \emph{repulsive tangent-point potential} 
\[
\TP[y] = \frac{2^{-q}}{q} \int_\o \int_\o \frac{1}{r^q(y(x),y(\tx))}\dv{\tx}\dv{x},
\]
where  $r(y(x),y(\tx))$ denotes the the radius of the sphere that is tangent to the deformed surface $y$ in the point $y(x)$ and which intersects $y$ in $y(\tx)$, cf. Figure~\ref{fig:radius}.
With $\nu_y(x) = \p_1 y(x) \times \p_2 y(x)$ the unit normal to the deformed surface at $y(x)$, this radius can be explicitly computed via 
\[
r(y(x),y(\tx)) = \frac{|y(x)-y(\tx)|^2}{2|\nu_y(x)\cdot(y(x)-y(\tx))|},
\]
which tends to zero whenever two distinct points on the deformed surface approach each other, thus causing a singularity in the potential for non-physical deformations with self-intersections. The use of the tangent-point potential was proposed for self-avoiding curves and surfaces in \cite{GoMa99,BaGoMadMar03} and successfully used in the simulation of knots \cite{BaReRi18}. 
It is further known that for choices of $q > 4$ the potential is self-avoiding for $2$-dimensional sub-manifolds of $\R^3$, i.\,e. we have $\TP[y] = \infty$ whenever $y$ is not one-to-one \cite{StMo11,Blatt13,Reiter17}. 
Consequently, after including the tangent-point potential in the energy functional and, thus, penalizing any tendency towards self-contact, we expect minimizers of finite energy to not show any self-intersections in a continuous setting.
Since the inverse  of $r(y(x),y(\tx))$ is an approximation of the normal curvature in the tangential direction $y(x)-y(\tx)$ if $|x - \tx| \ll 1$, and since we expect curvature to be bounded in some sensible way due to the nature of our model, we argue that we may exclude the singular values from the inner integral in the tangent-point potential and instead consider the simplified potential
\[
\TP_\veps[y] = \frac{2^{-q}}{q} \int_\o \int_{\o_\veps(x)} \frac{1}{r^q(y(x),y(\tx))}\dv{\tx}\dv{x},
\]
where $\o_\veps(x) = \o \setminus B_\veps(x)$ for some suitably chosen $\veps > 0$.
The self-avoiding bending energy is then defined for a parameter $\rho > 0$  via
\begin{equation*}\label{eq:energy}
E[y] = E_\mathrm{bend}[y] + \rho \TP_\veps[y],
\end{equation*}
where $E_\mathrm{bend} = E_\mathrm{Ki}$ or $E_\mathrm{bend} = E_\mathrm{bil}$, and we seek minimizers $y$ of $E[y]$ in the set of admissible functions
\[
\cA = \big\{ y \in H^2(\o; \R^3) : y|_\GD = y_\DD, \; \nabla y|_\GD = \phi_\DD,\; [\nabla y]\top \nabla y = \Id_2 \; \text{ a.\,e. in } \o \big\}.
\]
An intuitive way of thinking is that the parameter $\rho$ induces a characteristic length scale that defines a minimal positive distance. The density result from~\cite{Hornung08} implies that any function in $\cA$ can be approximated with arbitrary precision by smooth isometries, although an additional assumption is needed to guarantee compatibility of boundary conditions. The nonlinear isometry constraint in the definition of $\cA$ requires special attention in numerical approximations. Our approach employs discrete Kirchhoff triangles in the spatial discretization and then, in the single-layer case, uses the semi-implizit discrete gradient flow
\[ 
(d_t y^k,\vphi )_{H^2} + (D^2y^k,D^2\vphi) = (f,\vphi)-\rho\TP_\veps'[y^{k-1},\vphi] 
\]
to detect and approximate critical points of the resulting discrete energies, while enforcing a linearization of the isometry constraint in every step. In the bilayer case, following the approach which has been discussed in~\cite{BaPa20}, we use the property of isometries, $|\mathit{II(y)}| = |D^2y|$, to rewrite the energy as
\[
E_\mathrm{bil}[y] = \frac{1}{2} \int_\o |D^2 y|^2 \dv{x} - \a \int_\o \Delta y \cdot [\p_1 y \times \p_2 y] \dv{x} + a^2 |\o| - \int_\o f \cdot y \dv{x}
\]
and then, in the discretization of the resulting gradient flow, treat the first term, which is convex, implicitly, while treating the second, nonlinear term explicitly.

\section{Preliminaries}

\subsection{Variation of the tangent-point potential}\label{sec:variationTP}

The variational derivative of $\TP_\veps$ can be computed using chain, quotient and product rules and is given by
\begin{equation*}
\TP_\veps'[y;\vphi] = \int_\o \int_{\o_\veps(x)} \frac{|\nu_y \cdot (y(x)-y(\tx))|^{q-1}}{|y(x)-y(\tx)|^{2q+2}} \Bigl[ \cB[y,\vphi](x,\tx) + \cC[y,\vphi](x,\tx) - \cD[y,\vphi](x,\tx) \Bigr] \dv{\tx} \dv{x},
\end{equation*}
where
\begin{equation}\label{eq:tpvariation}
\begin{split}
&\cB[v,w](x,\tx) = \sign\Bigl(\nu_v(x) \cdot \bigl(v(x)-v(\tx)\bigr)\Bigr)\bigl|v(x)-v(\tx)\bigr|^2\\
&\hspace{2.5cm}\Bigl[ \bigl(\p_1 w(x) \times \p_2 v(x) + \p_1 v(x) \times \p_2 w(x) \bigr) \cdot \bigl(v(x)-v(\tx)\bigr) \Bigr],\\
&\cC[v,w](x,\tx) = \sign\Bigl(\nu_v(x) \cdot \bigl(v(x)-v(\tx)\bigr)\Bigr)\bigl|v(x)-v(\tx)\bigr|^2\\
&\hspace{2.5cm}\Bigl[ \nu_v(x) \cdot\bigl(w(x)-w(\tx)\bigr) \Bigr],\\
&\cD[v,w](x,\tx) = 2 \Bigl| \nu_v(x) \cdot \bigl(v(x)-v(\tx) \bigr) \Bigr| \Bigl( \bigl (v(x)-v(\tx) \bigr) \cdot \bigl(w(x)-w(\tx)\bigr) \Bigr).
\end{split}
\end{equation}

\subsection{Approximation spaces}

In the spatial discretization we avoid an $H^2$-conforming finite element method by employing discrete Kirchhoff triangular (DKT) elements. The discrete function spaces are subspaces $W_h \subset H^1(\o,\R)$ and $\Theta_h\subset H^1(\o,\R^2)$, corresponding to a triangulation $\cT_h$ of $\o$ into triangles with maximal diameter $h>0$.
For the triangulation $\cT_h$, we let $\cN_h$ and $\cS_h$ denote the set of vertices and sides of elements, respectively. The approximation spaces are then defined via
\begin{equation*}
\begin{split}
W_h =& \left\{ w_h \in C(\bar \o,\R) : w_h|_T\in P^{\mathrm{red}}_3(T) \text{ for all } T \in \cT_h ,\ \nabla w_h \text{ continuous in all }  z\in \cN_h \right\},\\
\Theta_h =& \left\{ \t_h\in C(\bar \o,\R): \t_h|_T\in P_2(T) \text{ for all } T\in \cT_h \right\},
\end{split}
\end{equation*}
where $P_k(T)$ denotes the set of polynomials of degree less or equal to $k\geq 0$ restricted to $T$ and $P^{\mathrm{red}}_3$ denotes the subset of cubic polynomials on $T$ defined by
\[ P^{\mathrm{red}}_3(T) = \left\{ p\in P_3(T) : p(x_T) = \sum_{z\in \cN_h\cap T} \left[ p(z) + \nabla p(z) \cdot (x_T - z) \right] \right\}, \]
with $x_T = (1/3) \sum_{z\in \cN_h \cap T} z$ the center of mass of T, i.\,e. $P^{\mathrm{red}}_3(T)$ results from $P_3(T)$ by eliminating one degree of freedom. A canonical interpolation operator $\cI_{W_h}:C^1(\bar \o) \to W_h$ for continuously differentiable functions is well defined via the identities $\cI_{W_h}w(z) = w(z)$ as well as $\nabla\cI_{W_h}w(z) = \nabla w(z)$ for all nodes $z\in \cN_h$.

The approximation of bending deformations with DKT elements is based on the construction of a discrete gradient operator
\[ \nabla_h : W_h \to \Theta_h^2 \]
which allows the definition of discrete second order derivatives of functions $w_h \in W_h$ via
\[ D^2_h w_h = \nabla\nabla_h w_h. \]

The degrees of freedom in $W_h$ are the function values and the derivatives at the vertices of the elements, whereas the degrees of freedom in $\T_h$ are the function values at both the vertices and the midpoints of element sides.

Let $S \in \cS_h$ be a side with endpoints $z_S^1, z_S^2$ and $z_S = \frac 1 2 (z_S^1+z_S^2)$ be its midpoint.
Denote $t_S$ a normalized tangent vector and $n_S$ a unit normal to the side $S$.
For $w_h\in W_h$ the discrete gradient operator is the uniquely defined, piecewise quadratic, continuous vector field $\theta_h\in \Theta_h^2$, such that for every node $z\in\cN_h$ and every side $S\in\cS_h$ the vector field satisfies the conditions
\begin{equation*}
  \begin{split}
    \theta_h(z) &= \nabla w_h(z),\\
    \theta_h(z_S)\cdot t_S &= \nabla w_h(z_S) \cdot t_S,\\
    \theta_h(z_S)\cdot n_S &= \frac 1 2 (\nabla w_h(z_S^1) + \nabla w_h(z_S^2)).
  \end{split}
\end{equation*}
\begin{remark}\label{rem:nablah}
  The discrete gradient operator satisfies the following approximation properties for all $w\in H^3(\o)$, $w_h\in W_h$ and $T\in \cT_h$, cf.~\cite[$\S 5$]{Braess07}:\\
  (i) There exists $c_1>0$ such that we have for $\ell=0,1$
  \[
  c_1^{-1} \|\nabla^{\ell+1}w_h\|_{L^2(T)} \leq \|\nabla^\ell \nabla_h w_h\|_{L^2(T)} \leq c_1 \|\nabla^{\ell+1} w_h\|_{L^2(T)},
  \]
  where $\nabla^1 = \nabla$ and $\nabla^0 = I$.\\
  (ii) There exists $c_2>0$ such that
  \[
  \|\nabla_h w - \nabla w\|_{L^2(T)} + h_T\|\nabla\nabla_h w - D^2 w\|_{L^2(T)} \leq c_2 h^2_T \|D^3 w\|_{L(T)}.
  \]
  (iii) There exists $c_3>0$ such that
  \[
  \|\nabla_h w_h - \nabla w \|_{L^2(T)} \leq c_3 h_T \|D^2 w_h\|_{L^2(T)}.
  \]
\end{remark}
\begin{remark}
  As a consequence of the last inequality in Remark~\ref{rem:nablah}, the mapping $w_h \mapsto \|\nabla\nabla_h w_h\|_{L^2(T)}$ defines a semi-norm on $W_h$, as well as a norm on every subspace of $W_h$ with vanishing values and derivatives on $\Gamma_D$.
\end{remark}

\section{Discretization and Minimization}

\subsection{Discrete energy}

We denote the nodal interpolation operator into continuous piecewise linear functions~$\cI_h^1$, the element-wise nodal interpolation operator into piecewise linear functions~$\widetilde \cI_h^1$ and the lumped $L^2$ inner product $(v,w)_h = \int_\o \widetilde \cI_h^1[v \cdot w]\dv{x}$. Assuming that the body force $f$ is piecewise continuous, the discrete energy is defined as
\begin{equation}\label{eq:discreteEnergy}
    E_h[y_h] =  E_h^{\mathrm{bend}}[y_h] + \rho \TP_{h}[y_h],
\end{equation}
with the elastic energy
\begin{equation*}
    E_h^{\mathrm{bend}}[y_h]  = \frac 1 2 \int_\o |\nabla\nabla_h y_h|^2 \dv{x} -(f,y_h)_h,
\end{equation*}
or 
\begin{equation*}
    E_h^{\mathrm{bend}}[y_h]  = \frac 1 2 \int_\o |\nabla\nabla_h y_h|^2 \dv{x} - \a\int_\o   \widetilde \cI_h^1 \bigl[ (\diver \nabla_h y_h) \cdot [\p_1 y_h \times \p_2 y_h] \bigr] \dv{x} + \a^2|\o|-(f,y_h)_h,
\end{equation*}
in the single and bilayer case, respectively, and the discrete self-avoidance functional
\begin{equation*}
    \TP_{h}[y_h] = \frac{2^{-q}}{q}\int_\o \cI_h^1 \bigg[ \int_{\o_{h}(x)} \widetilde \cI_h^1\bigg[ \frac 1 {r\bigl(y_h(x),y_h(\tx)\bigr)^q} \bigg] \dv{\tx}  \bigg] \dv{x},
\end{equation*}
where 
\begin{equation*}
  \o_{h}(x) = \o \setminus \bigcup \{ T \in \Th : B_{(h_\mathrm{min}/2)}(x) \cap T \ne \emptyset \}
\end{equation*}
is the union of all triangles which have an empty intersection with the $(h_\mathrm{min}/2)$ ball around $x$, with $h_\mathrm{min}$ denoting the smallest element diameter.
This choice guarantees that we have $\o_h(z) = \bigcup \{ T \in \Th : z \ne T \}$ for every vertex $z \in \Nh$.

Let $W_{h,D}$ denote the subset of $W_{h}$ satisfying given Dirichlet boundary conditions, i.e.
\[
W_{h,D} = \big\{w_h\in W_h(\cT_h) : w_h(z) =  w_D(z), \nabla w_h(z) = \phi_D(z) \text{ for all } z\in \cN_h \cap \Gamma_D \big\}
\]
and
\[
W_{h,0} = \big\{ w_h \in W_{h}(\cT_h) : w_h(z) = 0,\ \nabla w_h(z) = 0 \text{ for all } z \in \cN_h\cap \Gamma_D\big\}.
\]
The set $\cA_h$ of admissible discrete deformations is then defined via
\begin{equation*}
  \begin{split}
    \cA_h = \big\{ w_h \in W_{h,D}^3(\cT_h) : [\nabla w_h(z)]^\top \nabla w_h(z) = I_2 \text{ for all } z\in \cN_h \big\}
  \end{split}
\end{equation*}
and its tangent space at $y_h\in\cA_h$ is given by
 \begin{equation*}
  \begin{split}
    \cF_h[y_h] = \big\{ w_h \in W_{h,0}^3(\cT_h) : [\nabla w_h(z)]^\top \nabla y_h(z) + [\nabla y_h(z)]^\top \nabla w_h(z)  = 0 \text{ for all } z\in \cN_h, \big\}.
  \end{split}
\end{equation*}

\subsection{Discrete gradient flow}

For the minimization of the discrete energies~\eqref{eq:discreteEnergy} in $\cA_h$, we employ a semi-implicit discrete gradient flow scheme which is constrained to the linearization of the isometry constraint in every pseudo time step. 
As a consequence the iterates in the discrete gradient flow do not satisfy the discrete isometry constraint and numerically computed approximate minimizers are elements of the relaxed admissible set
\[
\cA_h^\d = \big\{ w_h \in W_{h,D}^3(\cT_h) : \big|[\nabla w_h(z)]^\top \nabla w_h(z) - I_2\big| \le \delta_h \text{ for all } z\in \cN_h \big\}
\]
that allows for some tolerance $\d_h > 0$, cf.~Remark~\ref{rem:tol}.
Practical properties of the proposed method are investigated in the numerical experiments in Section~\ref{sec:experiments}.

We use the notation $(\cdot\,,\cdot)_\ast = (\nabla\nabla_h\,\cdot\,,\nabla\nabla_h\,\cdot)$ and $\|\cdot\|_\ast$ to denote the discrete $H^2$ scalar product which we use to define the gradient flow and its induced norm, respectively.
With the functionals $\cB$, $\cC$ and $\cD$, defined as in~\eqref{eq:tpvariation}, the variational derivative $\TP_{h}'[y_h; \vphi_h]$ of the discrete tangent-point potential $\TP_{h}$ is given by
\begin{equation}\label{eq:tpvarapprox}
\TP_{h}'[y_h; \vphi_h] = \int_\o \cI_h^1 \cQ_h[y_h,\vphi_h] \dv{x},
\end{equation}
with
\[
\cQ_h[y_h, \vphi_h](x) = \int\limits_{\o_{h}(x)} \widetilde \cI_h^1\Bigl[ R\big[y_h\big]\big(x,\tx\big) \Big( \cB[y_h,\vphi_h](x,\tx) +\cC[y_h,\vphi_h](x,\tx) - \cD[y_h,\vphi_h](x,\tx) \Big) \Bigr]\dv{\tx},
\]
where
\[
R\big[y\big]\big(x,\tx\big) = \frac{\bigl|\nu_{y}(x)\cdot \big(y(x)-y(\tx)\big)\bigr|^{q-1}}{\bigl|y(x)-y(\tx)\bigr|^{2q+2}}.
\]
We can now formulate the linearly constrained discrete gradient flow for the minimization of~\eqref{eq:discreteEnergy}.

\begin{algorithm}[discrete isometry flow]\label{alg:isoflow}
  Given an initial deformation $y_h^0 \in \cA_h$, choose a step size $\tau>0$ and a stopping criterion $\veps_{\mathrm{stop}}>0$ and set $k=1$.\\
  (1) Compute $d_t y_h^k \in \cF_h[y_h^{k-1}]$ such that
    \begin{equation*}
        (d_t y_h^k,\vphi_h)_\ast = - \left( \nabla\nabla_h (y_h^{k-1} + \tau d_t y_h^k),\nabla\nabla_h\vphi_h \right) + \left( f,\vphi_h \right)_h + \rho \TP_{h}'[y_h^{k-1}; \vphi_h]
    \end{equation*}
    for all $\vphi_h \in \cF_h[y_h^{k-1}]$.\\
  (2) Set $y_h^k = y_h^{k-1}+\tau d_ty_h^k$.
    If $\|d_t y_h^k\|_\ast<\veps_{\mathrm{stop}}$, stop the iteration. Otherwise, increase $k$ via $k \mapsto k+1$ and continue with (1).
\end{algorithm}

\begin{remark}\label{rem:tol}
  Algorithm~\ref{alg:isoflow} does not contain any projection step which would guarantee the satisfaction of the discrete isometry constraint in the nodes of the triangulation. Omitting such a step is motivated by corresponding rigorous results on discrete gradient flows for single layer plates~\cite{Bartels13}, bilayer plates~\cite{BaBoNo17, BaPa20} as well as harmonic maps~\cite{Bartels16}. The violation of the respective constraints is independent of the number of performed iterations and controlled by the step size~$\tau$.
\end{remark}

\subsection{Implementation}
The linearization of the isometry constraint in combination with the explicit treatment of nonlinear parts of the energy functional in the discrete gradient flow lead to linear systems in every step. 
The use of the DKT element serves as a model discretization of the problem under consideration. The transfer of our concepts to more standard elements such as dG is straightforward.
If one has, however, obtained an implementation of the DKT element and its discrete gradient operator (see Section 8.2 of~\cite{Bartels15} for details on how to implement DKT), the implementation of Algorithm~\ref{alg:isoflow} is uncomplicated.  Denoting with $K = \# \Nh$ the number of vertices of the triangulation and using the nodal basis $[\psi_i]_{1 \le i \le 3K}$ of the discrete function space $W_h$, we can identify functions $y_h \in W_h^3$ with vectors $\mathbf{y} \in \R^{9K}$. The linearized discrete isometry constraint can then be imposed by introducing Lagrange multipliers, i.\,e. by considering the equivalent saddle point problem 
\begin{equation}\label{eq:saddle}
  \begin{bmatrix}
    (1+\tau) \mathbf{S}_\mathrm{DKT} & [\mathbf{A}^{k-1}]^\top \\
    \mathbf{A}^{k-1} & 0
  \end{bmatrix}
  \begin{bmatrix}
    d_t \mathbf{y}^k \\
    \mathbf{\lambda}
  \end{bmatrix}
  =
  \begin{bmatrix}
    - \mathbf{S}_\mathrm{DKT} \mathbf{y}^{k-1} + \mathbf{b}_f + \rho \mathbf{b}_{\TP'}^{k-1} \\
    0
  \end{bmatrix}
\end{equation}
in every time step. Here, the matrix $\mathbf{S}_\mathrm{DKT} \in \R^{9K \times 9K}$ encodes the scalar product $(\nabla \nabla_h \cdot,\nabla \nabla_h \cdot)$, the matrix $\mathbf{A}^{k-1} \in \R^{9K \times 3K}$ encodes the constraint map 
\[
w_h \mapsto \cI_h^1 \left[(\nabla w_h)^\top \nabla y_h^{k-1} + (\nabla y_h^{k-1})^\top \nabla w_h \right]
\]
and the two vectors 
\begin{equation*}
  \mathbf{b}_f = \begin{bmatrix} 
  (f, \psi_1)_h \\ \vdots \\ (f, \psi_{9K})_h
  \end{bmatrix}, \quad
  \mathbf{b}_{\TP'}^{k-1} = \begin{bmatrix} 
  \TP'_{h}[y_h^{k-1};\psi_1] \\ \vdots \\ \TP'_{h}[y_h^{k-1};\psi_{9K}]
  \end{bmatrix}
\end{equation*}
contain the contributions of the body force potential as well as the explicitly treated variation of the tangent-point potential, respectively. Note that the entries of the matrix $\mathbf{A}^{k-1}$, as well as the normals $\p_1 y_h(z) \times \p_2 y_h(z)$ needed in the computation of $\mathbf{b}_{\TP'}^{k-1}$, are directly obtained from the degrees of freedom in the employed DKT finite element space $W_h$.

A different, more elaborate, approach for treating the linear constraints involves the construction of a basis of the subspace $\ker \mathbf{A}^{k-1} \subset \R^{9K}$, such that the change of basis, $\mathbf{C} \colon \R^{6N} \to \ker \mathbf{A}^{k-1}$, leads to the symmetric positive definite system 
\begin{equation}\label{eq:red}
(1 + \tau) \mathbf{C}^\top \mathbf{S}_\mathrm{DKT} \mathbf{C} \hat{\mathbf{x}} = \mathbf{C}^\top \mathbf{b} 
\end{equation}
where, as before, the right-hand side is given by $\mathbf{b} = - \mathbf{S}_\mathrm{DKT} \mathbf{y}^{k-1} + \mathbf{b}_f + \rho \mathbf{b}_{\TP'}^{k-1}$. The required basis of $\ker \mathbf{A}^{k-1}$ can be obtained directly from the available degrees of freedom via a point-wise construction. This procedure leads to smaller linear systems with only half the number of unknowns compared to the approach based on Lagrange multipliers. Furthermore, the system matrix in~\eqref{eq:red} is s.p.d., whereas the system matrix in~\eqref{eq:saddle} does not have this property. Thus, problem~\eqref{eq:red} can possibly be addressed using an iterative solver, although the choice of an adequate preconditioner remains an open problem.

Considering that the computational cost in every time step is caused almost exclusively by the assembly of the variation of the tangent-point potential $\mathbf{b}_{\TP'}^{k-1}$, which has a complexity that grows at least quadratically in the number of elements, the gain in computing time resulting from smaller linear systems is negligible.

\section{Simplified assembly}

The assembly of the right-hand-side vector in the linear systems that arise in every time step of Algorithm~\ref{alg:isoflow} is computationally expensive. Element-wise quadrature rules for approximating the double integral in the tangent-point potential result in long summations with the number of summands growing quadratically with respect to the number of elements in a triangulation. Taking into account that specific problems might require a small mesh size to accurately resolve the problems' geometry and/or -- to prevent instabilities -- a relatively small (pseudo-)time step size for which hundreds of thousands of iterations are needed to reach a numerical equilibrium state, computation times may quickly become unacceptable.

In this section we propose, without further investigation, three different strategies which can be adopted individually or in arbitrary combinations to at least partially overcome the difficulties.

\subsection{Parallel computation}
\begin{figure}
  \includegraphics[height=5.7cm]{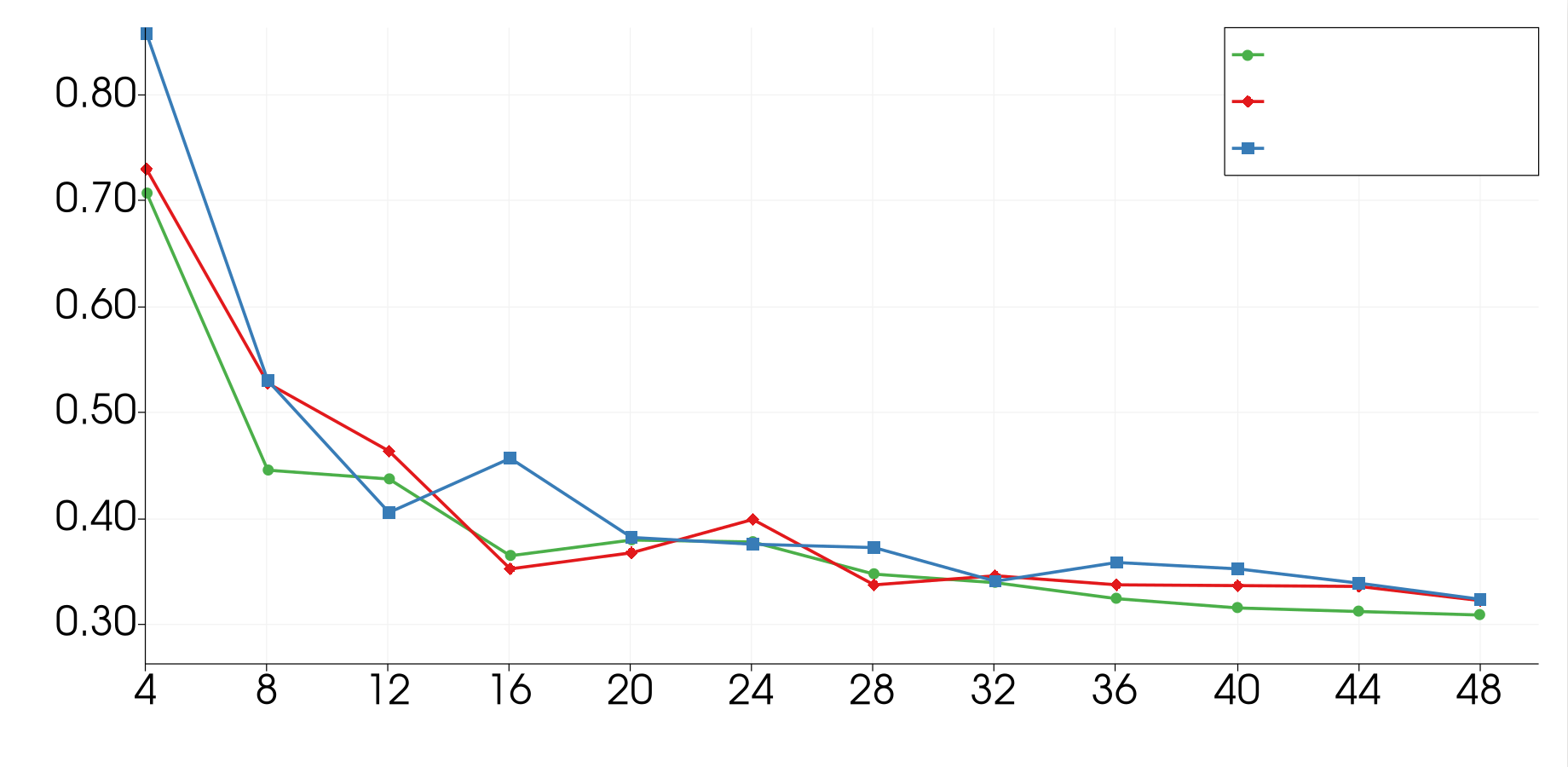}
  \put(-333,51){\rotatebox{90}{\small Parallel efficiency}}
  \put(-195,3){\small Number of processes}
  \put(-55,148){\footnotesize 1'600 trian.}
  \put(-55,138){\footnotesize 6'400 trian.}
  \put(-59,128){\footnotesize 25'600 trian.}
  \caption{Parallel efficiency measured for the assembly time of the linear systems in each pseudo time step of Algorithm~\ref{alg:isoflow} in Example~\ref{ex:trefoil} after parallelizing the computation of the right-hand-side vector $\mathbf{b}_{\TP'}^{k-1}$.}
  \label{fig:pareff}
\end{figure}
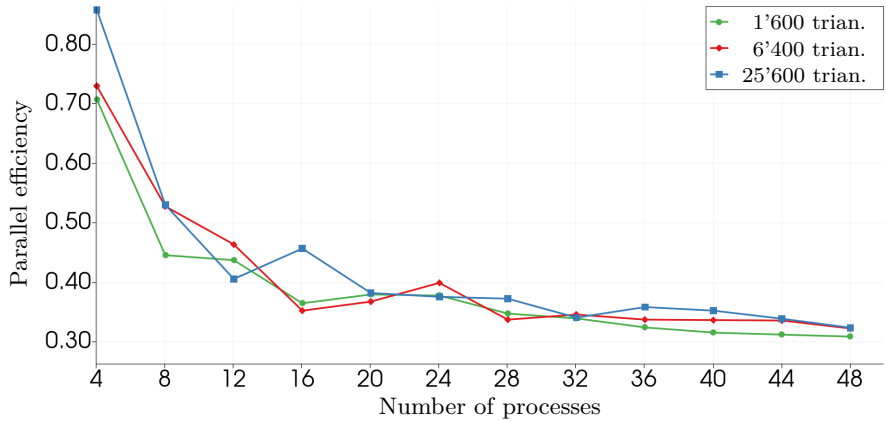

The components of the computationally expensive right-hand-side vector $\mathbf{b}_{\TP'}^{k-1}$ are given by integrals that can be naturally decomposed into subintegrals via
\[
\TP'_{h}[y_h^{k-1};\psi_i] = \sum_{T \in \Th} \int_T \cI_h^1 \cQ_h[y_h,\psi_i] \dv{x}.
\] 
Being mutually independent quantities, the values of the subintegrals can be computed in parallel with very little overhead, resulting in a significant speed up of the assembly routine for the linear system in every time step of Algorithm~\ref{alg:isoflow}. In the numerical experiments in Section~\ref{sec:experiments} we use OpenMP for a simple parallelization in the computation of $\mathbf{b}_{\TP'}^{k-1}$. The resulting parallel efficiency $\frac{\text{speed-up}}{\#\text{threads}}$ measured on our $24\;\times\;$Intel\textsuperscript{\tiny\textregistered} Xeon\textsuperscript{\tiny\textregistered} CPU E5-2695 v2 @ 2.40GHz machine for the assembly of the complete linear system in Example~\ref{ex:trefoil} is illustrated in Figure~\ref{fig:pareff} and supports the hypothesis that even basic parallelization techniques significantly decrease computation times.

\subsection{Hierarchical quadrature}
A natural property observed in deformations $y \colon \o \to \R^3$ with a tendency towards self-contact is localization of relevant contributions expressed via the existence of a subdomain $\widetilde{\o} \subset \o$ such that $|\widetilde{\o}| \ll |\o|$ and $\TP[y|_{\widetilde{\o}}] \approx \TP[y]$.
Using a hierarchical $N$-level approach we aim to identify such \emph{high-potential regions} of the domain in every step of Algorithm~\ref{alg:isoflow} and subsequently construct a non-conforming, locally refined triangulation $\cT_N$ which itself is a coarsening of $\Th$. The variation of the tangent-point potential is then approximated via quadrature on the elements of the obtained coarsened triangulation. Such an approach may significantly reduce the number of elements involved and, thus, computation time spent in the assembly for the right-hand-side vector.
\begin{algorithm}[Hierarchical quadrature]
    Choose a sequence $\cT_{h_1},\cT_{h_2},\dotsc,\cT_{h_N}=\cT_h$ of triangulations, such that $\cT_{h_{i+1}}$ is the result of a red refinement of $\cT_{h_i}$. Fix the parameter $\sigma\in (0,1]$,  set $i=1$ and $\cT_1 = \cT_{h_1}$.\\
    (1) Choose a subset $\widetilde{\cT}_{i} \subset \cT_{i}$, such that 
    \begin{align*}
      \sum_{T\in\widetilde\cT_i}  \int_T \cI_h^1 \bigg[ \int_{\o_{h_i}(x)} \widetilde \cI_{h_i}^1\bigg[ & \frac 1 {r\bigl(y_h(x),y_h(\tx)\bigr)^q} \bigg] \dv{\tx} \bigg ] \dv{x}  \\
      &\geq \hspace{8pt} \sigma \int_\o \cI_h^1 \bigg [ \int_{\o_{h_i}(x)} \widetilde \cI_{h_i}^1\bigg[ \frac 1 {r\bigl(y_h(x),y_h(\tx)\bigr)^q} \bigg] \dv{\tx}  \bigg ] \dv{x}, 
    \end{align*}
    (2) Define $\cT_{i+1} = \widetilde{\cT}^\mathrm{rfd}_{i} \cup \cT_i \setminus \widetilde{\cT}_i$, where the set $\widetilde{\cT}^\mathrm{rfd}_{i}$ contains exactly the red-refined elements from the subset $\widetilde{\cT}_i$ that was chosen in (1). \\
    (3) If $i+1 = N$, stop the algorithm and use (element-wise) quadrature on $\cT_N$.
    Otherwise, increase $i$ via $i\mapsto i+1$ and continue with (1). 
\end{algorithm}

\subsection{Boundary-domain potential}
For every  deformation $y$ that includes self-intersections there are distinct intersection points $x, x' \in \o$ for which we have that $y(x) = y(x')$. An observation that can be made for isometric deformations $\R^2 \supset \o \to \R^3$ is the fact that in many relevant examples with self-intersections we can identify at least one intersection point that lies on the boundary~$\p \o$ of the domain.
This motivates the use of a \emph{boundary-domain tangent-point potential},
\[
\TP_{\p \o, \veps}[y] = \frac{2^{-q}}{q} \int_{\o} \int_{\p_\veps \o (x)} \frac{1}{r^q(y(x),y(\tx))}\dv{\tx}\dv{\sigma(x)},
\]
where $\p_\veps \o (x) = (\p \o) \setminus B_\veps(x)$, and which acts as a repulsive potential on the domain versus its boundary, and which can be discretized analogously to $\TP_\veps [y]$. Using this approach the computational cost can be reduced, since quadrature for the outer (now one-dimensional) integral only requires the use of significantly less quadrature points compared to the two-dimensional integral in $\TP_\veps [y]$. Unless the examples were artificially constructed to exhibit self-intersections without intersection points on $\p \o$, we found that utilization of the boundary-domain potential lead to similar results as the full tangent-point potential.

\section{Numerical Experiments}
\label{sec:experiments}

The numerical experiments reported in this section investigate the practical properties of our algorithm and illustrate situations that benefit from its utilization as well as situations where difficulties are encountered. The code for the numerical experiments was written in C\texttt{++}, incorporating several routines provided by DUNE \cite{DUNEa, DUNEb, alugrid} and using a direct solver provided by \mbox{UMFPACK} \cite{umfpack} for the solution of the linear systems in every step. Visualizations of typical evolutions in the discrete gradient flow can be found online~\cite{BaPaAnim}. Our code was parallelized as described in Section~4.1. However, in order not to introduce additional approximation errors, we refrain from employing hierarchical quadrature, cf. Section~4.2. The boundary-domain potential, cf. Section~4.3, is only employed in special settings: a comparison of results obtained with the full tangent-point potential versus the boundary-domain potential is provided in the context of Example~\ref{ex:twist180}. The open-source application ParaView \cite{paraview} was used to visualize the computed discrete deformations. We note that all plots show $P_1$-interpolants of the respective discrete functions in the space $W_h$, i.\,e. we neglect degrees of freedom corresponding to the derivatives of the deformation at the nodes of a triangulation. In the following we denote with
\[
\d_\mathrm{iso}[y_h] = \| \cI_h^1[\nabla y_h^\top \nabla y_h] - I_2\|_{L^\infty}
\]
the isometry error of a discrete function and with $y_h^\infty$ the numerically computed equilibrium state of a given problem that we obtain with the stopping criterion $\veps_\mathrm{stop} = 10^{-3}$.
Whenever the self-avoidance potential was considered, the coloring of the deformed surfaces corresponds to the magnitude of the tangent-point potential density, i.\,e. the value of
\[
\tp_y(x) = \frac{2^{-q}}{q} \cI_h^1 \bigg[ \int_{\o_{h}(x)} \widetilde \cI_h^1\bigg[ \frac 1 {r\bigl(y_h(x),y_h(\tx)\bigr)^q} \bigg] \dv{\tx} \bigg].
\]
Furthermore, we define a piecewise linear vector field $f_\tp \colon \o \to \R^3$ via
\[
f_\tp(z) = \bigg[ \int_{\o_h(z)}  \widetilde \cI_h^1 \Big[ \sign \big( y_i(z) - y_i(\tx)\big) \frac{(\nu_y(z))_i(y_i(z) - y_i(\tx))}{|y(z)-y(\tx)|^{2}} \Big] \dv{\tx} \bigg]_{i=1,2,3}
\]
for all $z \in \Th$, which we use to visualize the pseudo force that is induced by the tangent-point potential. Note that this vector field is used for the sake of visualization only and does neither consider the exponent $q$ nor the scaling factor $\rho$.
Since we do not know any analytical minimizers, we use the quantity
\[
\mathrm{EOC} = \log_2 \Bigg( 
  \frac{\|\cI_h^1 y_h^\infty - \cI_{h/2}^1 y_{h/2}^\infty\|_{L^2}}
  {\|\cI_{h/2}^1 y_{h/2}^\infty - \cI_{h/4}^1 y_{h/4}^\infty\|_{L^2}}
\Bigg)
\]
as an experimental order of convergence, computed from the final iterates $y_h^\infty$ obtained on the three highest refinement levels in each corresponding experiment.

\subsection{Simple compression of a strip}
We consider compressive boundary conditions on the short ends of a rectangular $10 \times 1$ single layer plate. 
The compressive boundary conditions are chosen such that a self-intersection occurs in numerical approximations of stationary points in the absence of a self-avoidance potential and, thus, that self-contact of the deformed deformed plate can be assumed to occur in a physical solution.
\begin{figure}
  \subfloat{
    \begin{overpic}[trim = 40mm -15mm 40mm -10mm, clip, height=40mm]{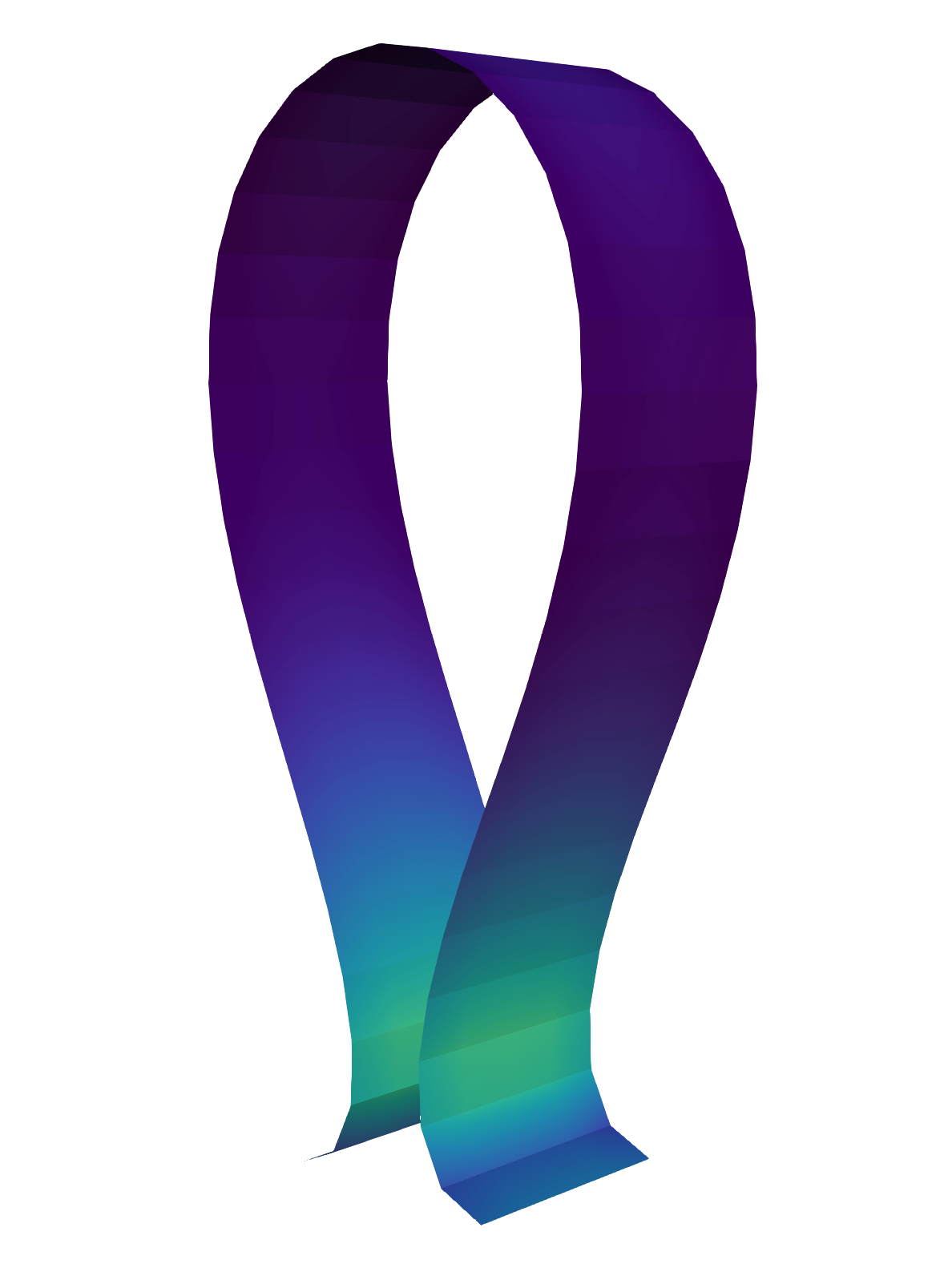}
      \put(5,0){\small $\#\cT_h = 320$}
      \put(5,-12){\small $q=5,\enspace\beta=1.0$}
    \end{overpic}
  }
  \subfloat{
    \begin{overpic}[trim = 40mm -15mm 40mm -10mm, clip, height=40mm]{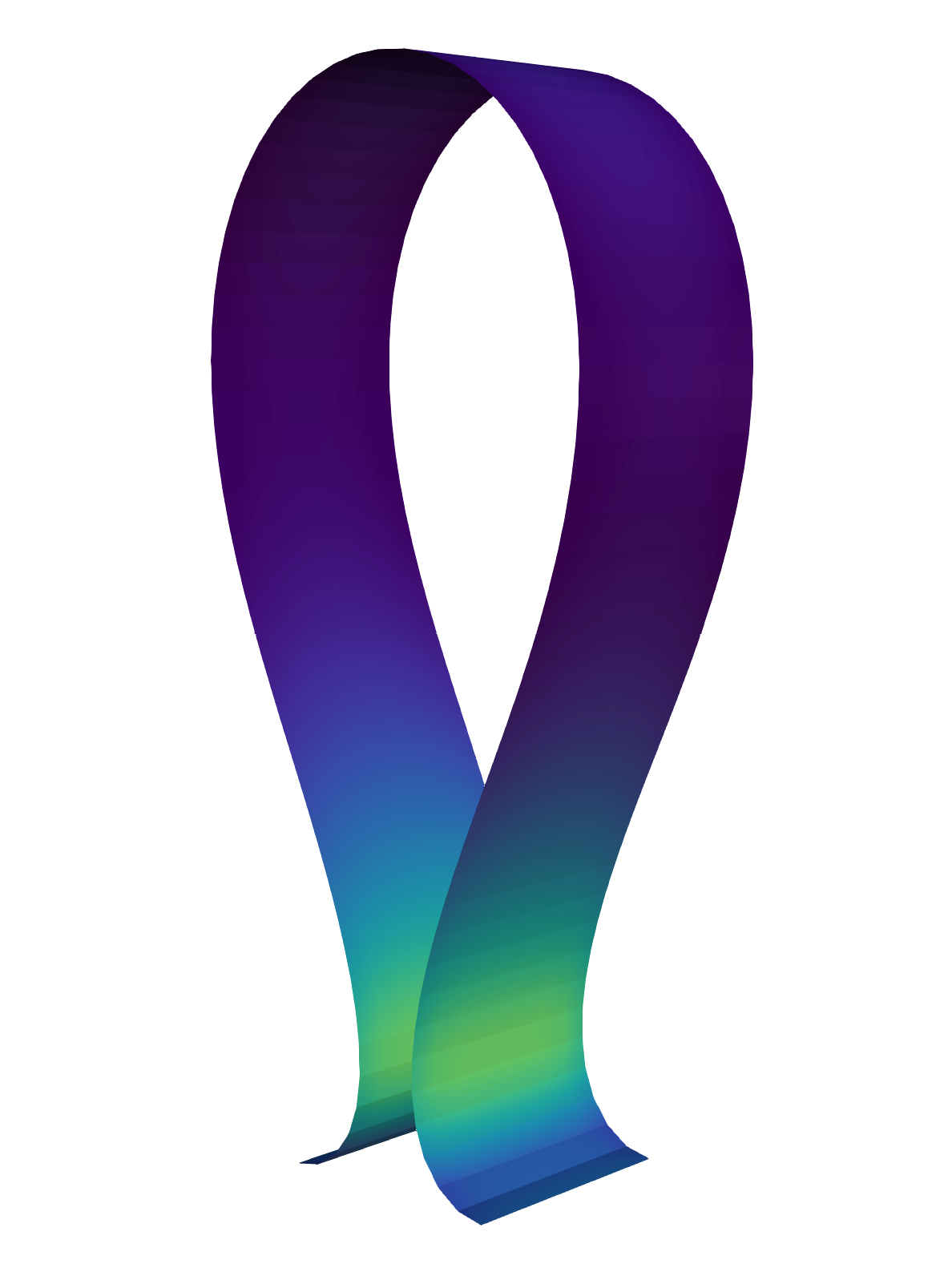}
      \put(5,0){\small $\#\cT_h=1280$}
      \put(5,-12){\small $q=5,\enspace\beta=1.0$}
    \end{overpic}
  }
  \subfloat{
    \begin{overpic}[trim = 40mm -15mm 40mm -10mm, clip, height=40mm]{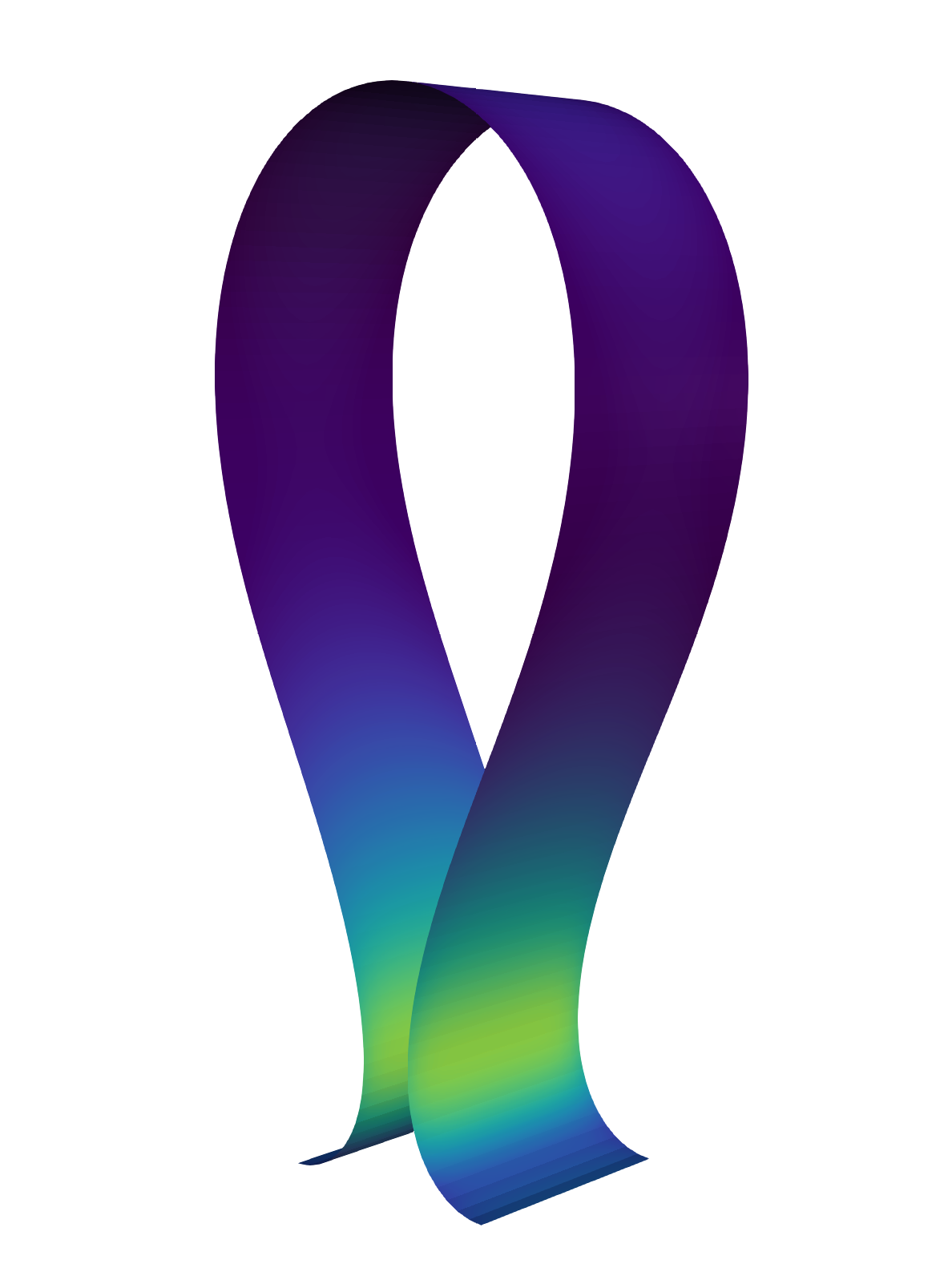}
      \put(5,0){\small $\#\cT_h=5120$}
      \put(5,-12){\small $q=5,\enspace\beta=1.0$}
    \end{overpic}
  }
  \subfloat{
    \begin{overpic}[trim = 40mm -15mm 40mm -10mm, clip, height=40mm]{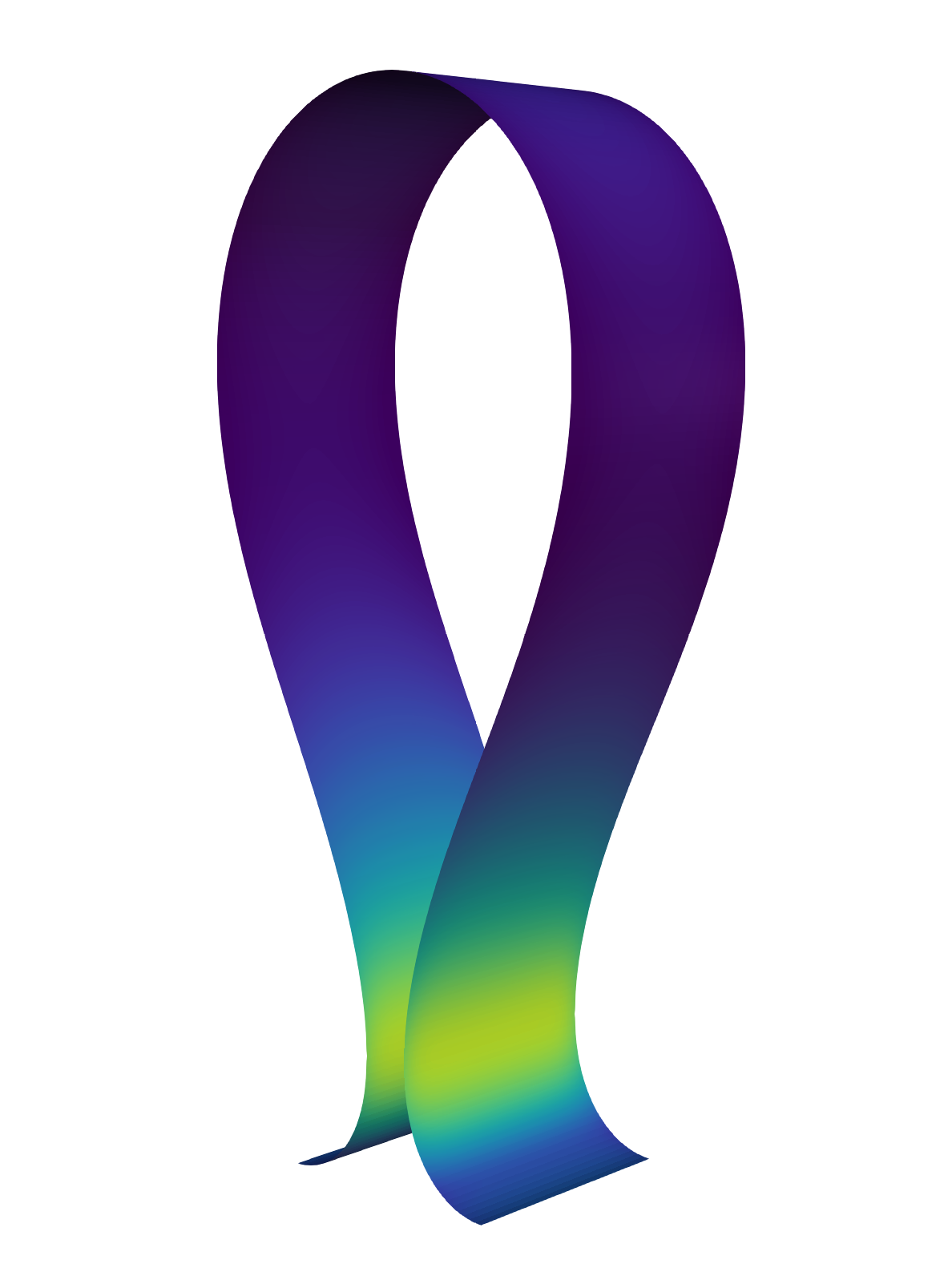}
      \put(5,0){\small $\#\cT_h=20480$}
      \put(5,-12){\small $q=5,\enspace\beta=1.0$}
    \end{overpic}
  }
  \\
  % NEWLINE
  \subfloat{
    \begin{overpic}[trim = 40mm -15mm 40mm -10mm, clip, height=40mm]{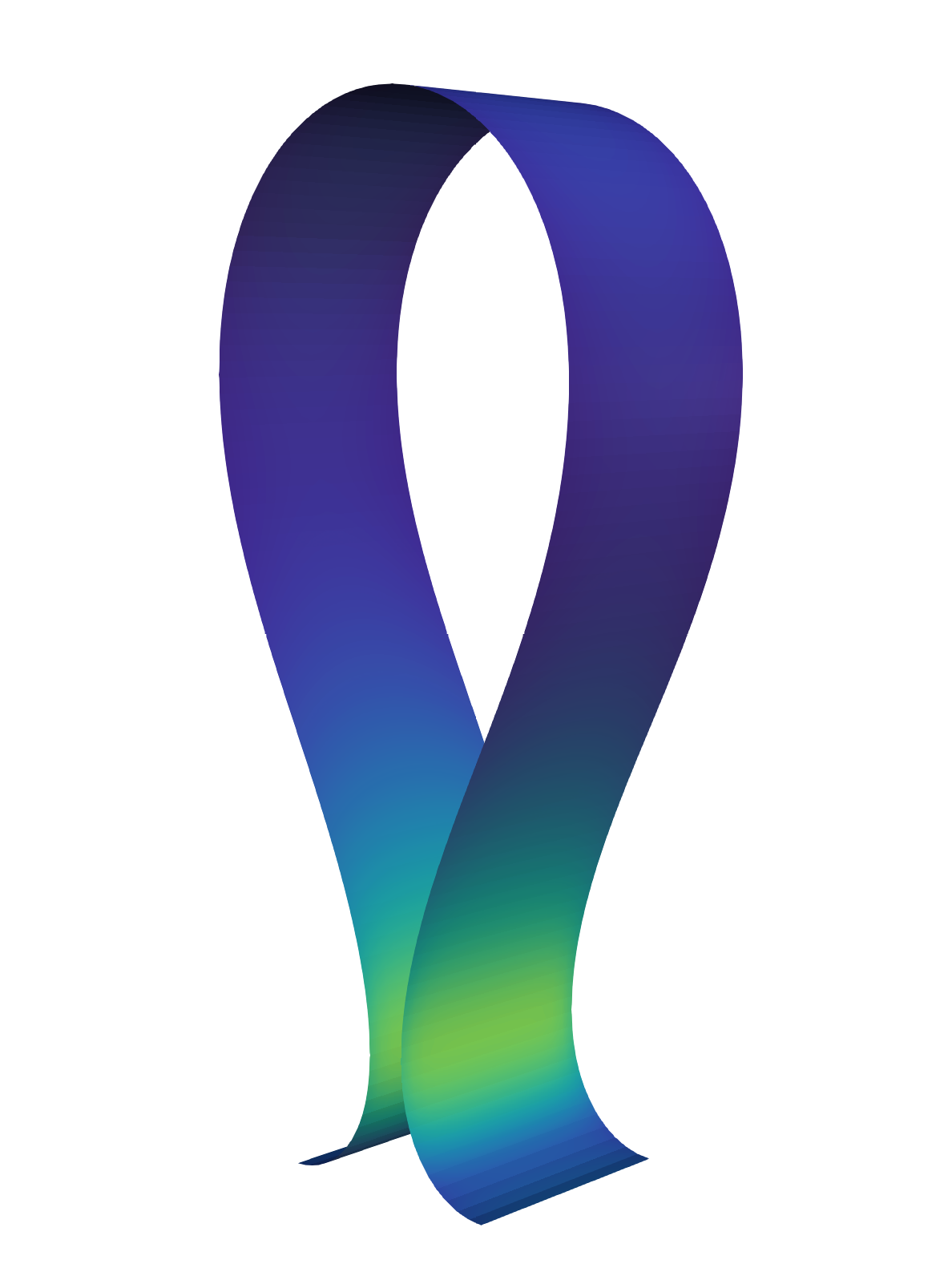}
      \put(5,0){\small $\#\cT_h=5120$}
      \put(5,-12){\small $q=4,\enspace\beta=1.0$}
    \end{overpic}}
  \subfloat{
    \begin{overpic}[trim = 40mm -15mm 40mm -10mm, clip, height=40mm]{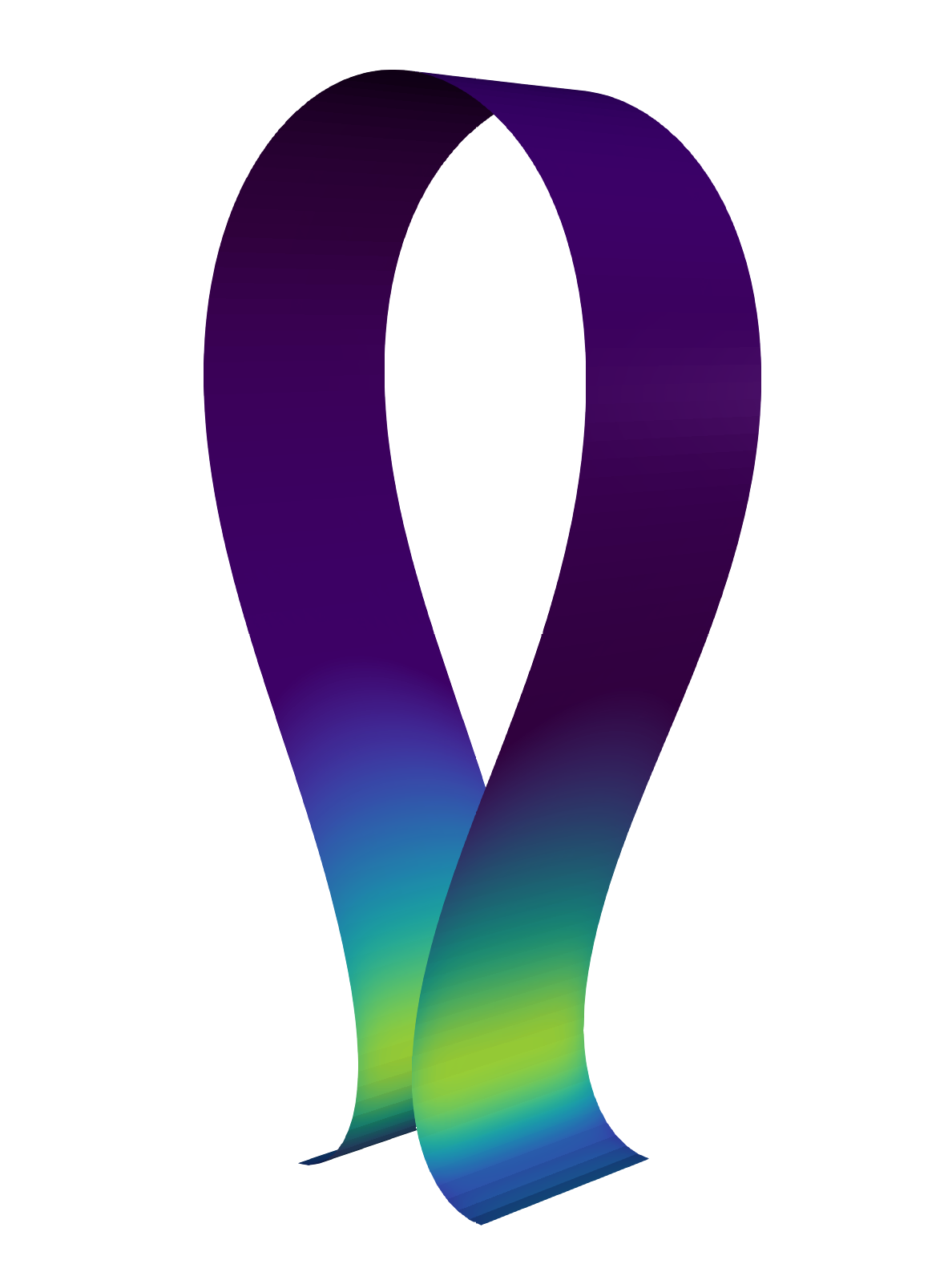}
      \put(5,0){\small $\#\cT_h=5120$}
      \put(5,-12){\small $q=6,\enspace\beta=1.0$}
    \end{overpic}
  }
  \subfloat{
    \begin{overpic}[trim = 40mm -15mm 40mm -10mm, clip, height=40mm]{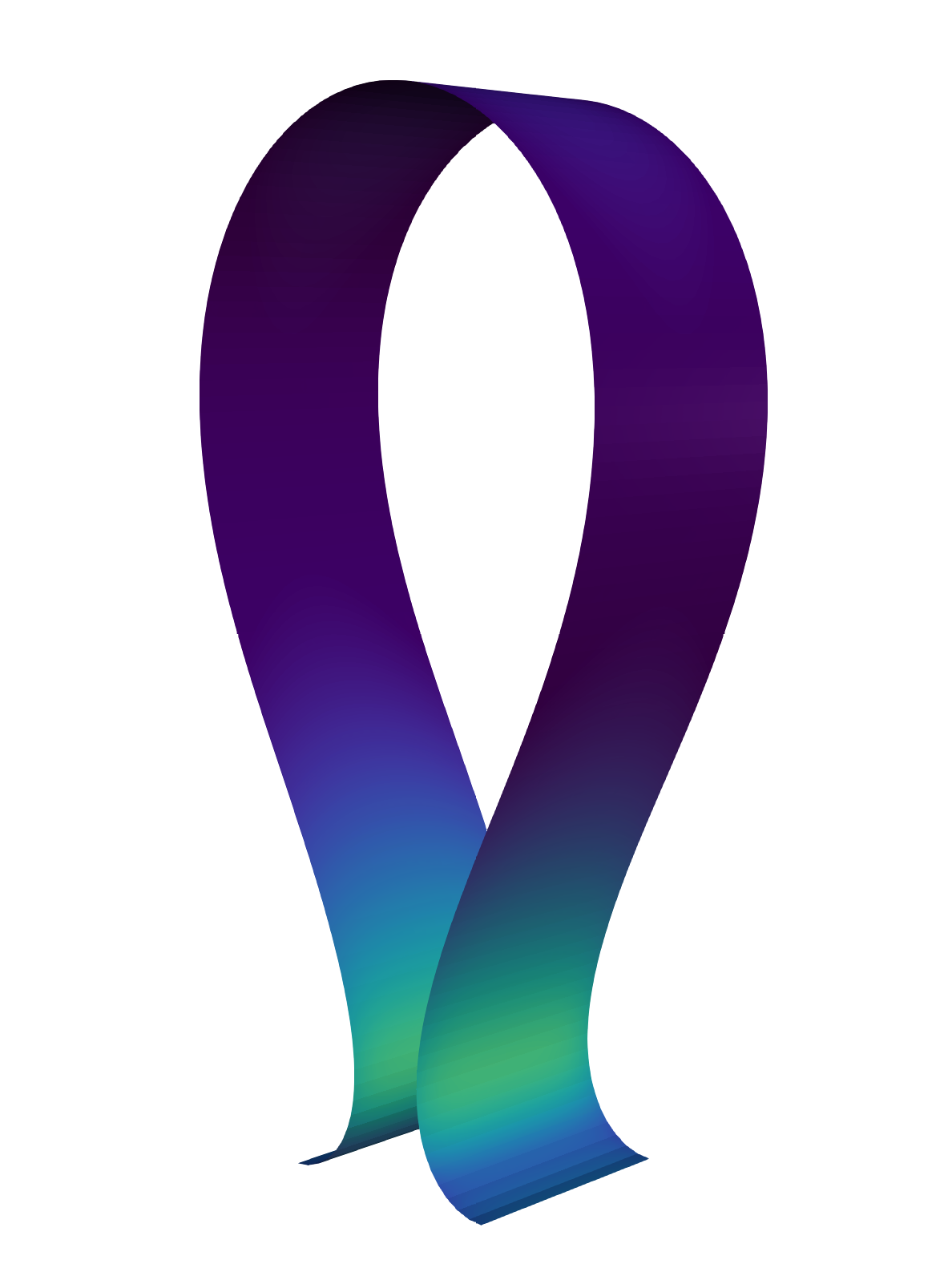}
      \put(5,0){\small $\#\cT_h=5120$}
      \put(5,-12){\small $q=5,\enspace\beta = 0.5$}
    \end{overpic}
  }
  \subfloat{
    \begin{overpic}[trim = 40mm -15mm 40mm -10mm, clip, height=40mm]{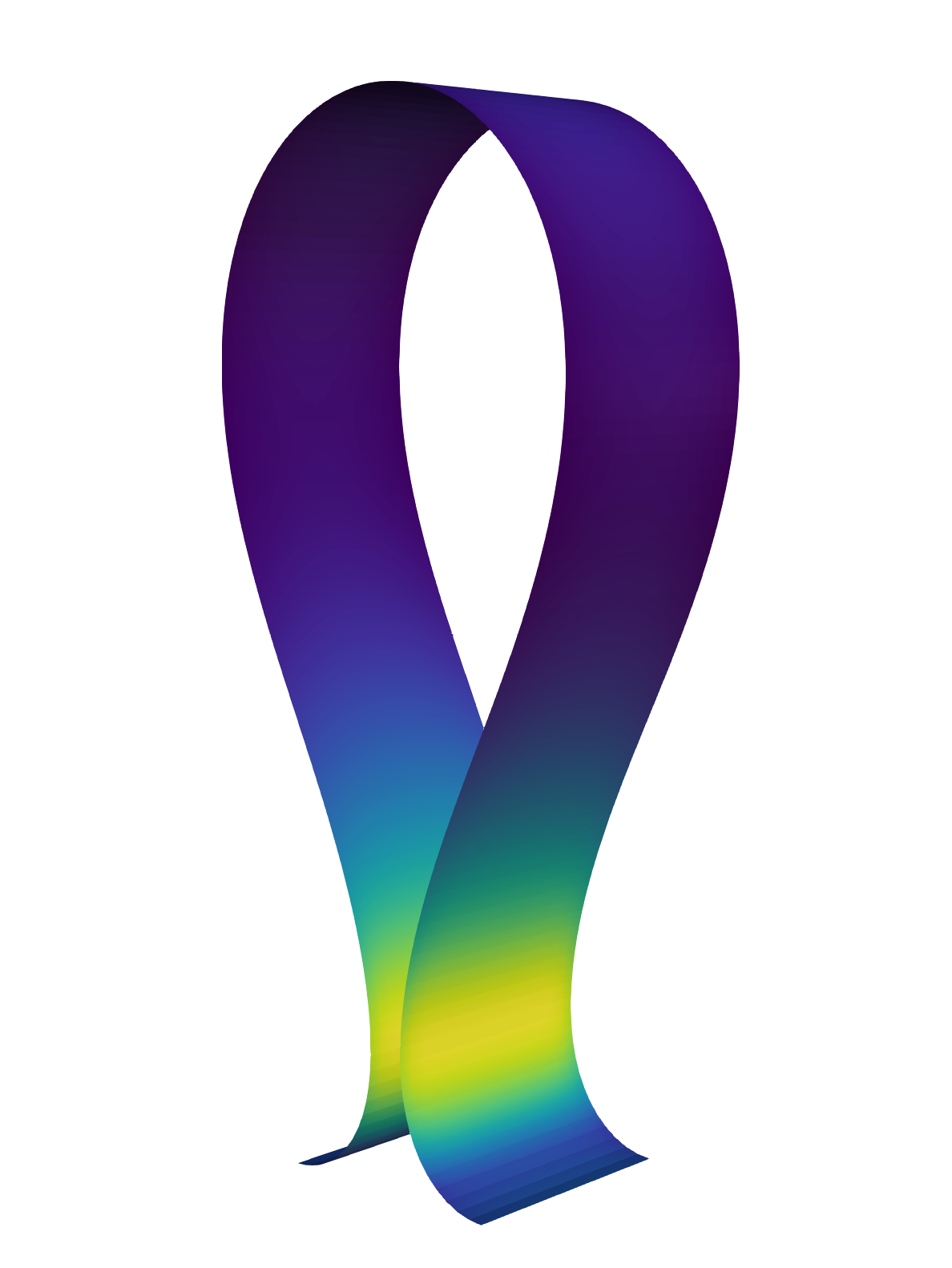}
      \put(5,0){\small $\#\cT_h=5120$}
      \put(5,-12){\small $q=5,\enspace\beta = 1.5$}
    \end{overpic}
  }
  \subfloat{
    \begin{overpic}[trim = 40mm 11mm 20mm 1mm, clip, height=40mm]{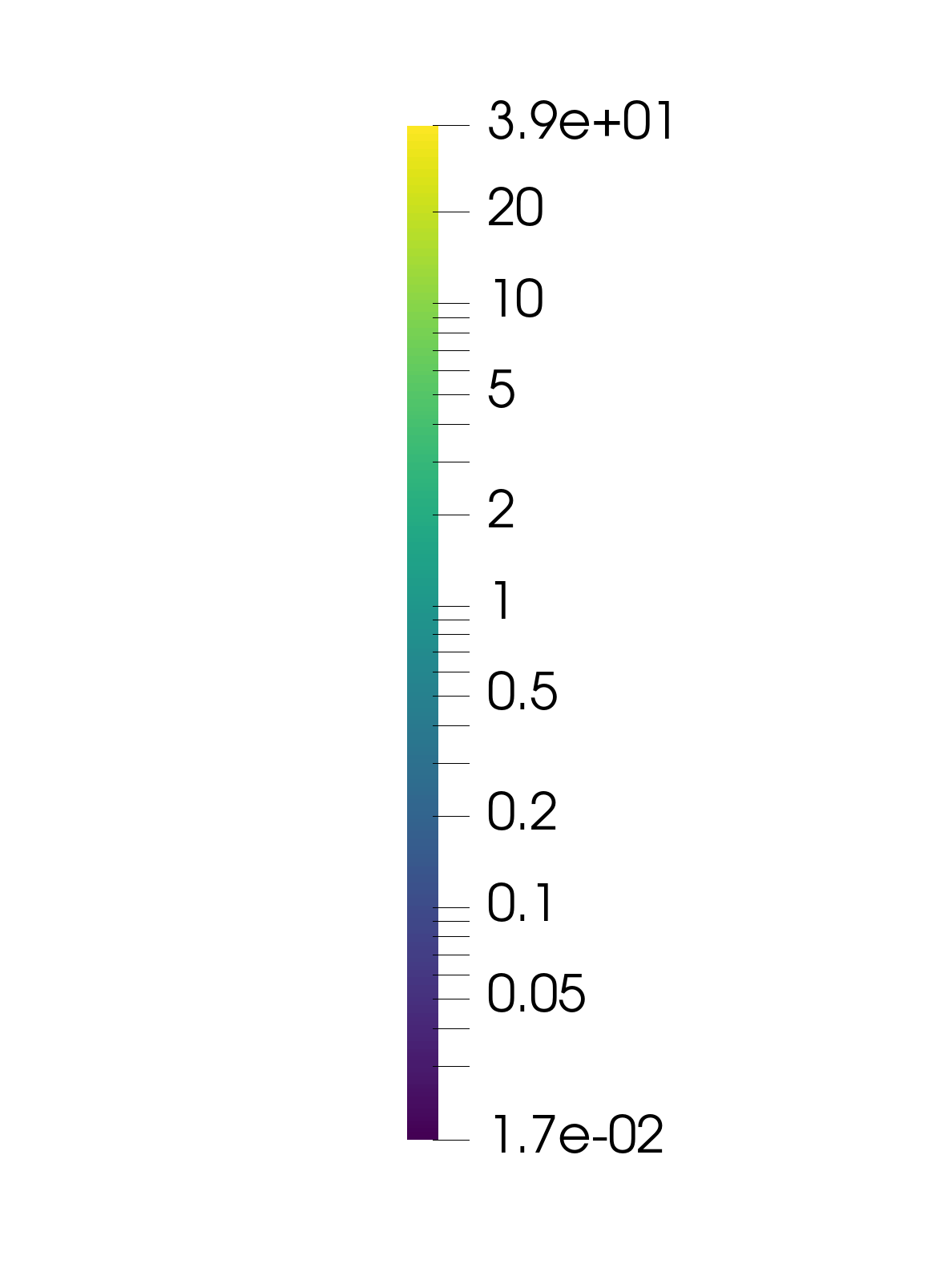}
      \put(27,56){\footnotesize $\mathrm{tp}_y(x)$}
    \end{overpic}
  }
  \\
  \caption{Stationary configurations of the self-avoiding compressed strip in Example~\ref{ex:compress} for $q=5$, $\beta=1$ on different triangulations and for different values of parameters $\rho = (\hat{h}/2)^\b$ and $q$ on the triangulation~$\mathcal{T}_4$.}
  \label{fig:compress}
\end{figure}
\begin{table}
  \centering
  \begin{tabular}{c c c c c c c}\hline
  $k$ & $\b$ & $q$ & $N_\mathrm{iter}$ & $E_h[y_h^\infty]$ & $\TP_{h}[y_h^\infty]$ & $\delta_\mathrm{iso}[y_h^\infty]$ \\ \hline
  $2$ & 1.0 & 5 & 448  & \num{6.61648} & \num{4.14697} & \num{0.171152} \\
  $3$ & 1.0 & 5 & 922  & \num{6.47112} & \num{7.18907} & \num{0.1498} \\
  $4$ & 1.0 & 4 & 1932 & \num{5.91176} & \num{10.4353} & \num{0.0570355} \\
  $4$ & 0.5 & 5 & 2667 & \num{7.20374} & \num{4.03493} & \num{0.0835602} \\
  $4$ & 1.0 & 5 & 2068 & \num{6.19095} & \num{11.1764} & \num{0.0573204} \\
  $4$ & 1.5 & 5 & 1932 & \num{5.73538} & \num{33.4938} & \num{0.0570295} \\
  $4$ & 1.0 & 6 & 2884 & \num{6.72736} & \num{12.2858} & \num{0.12109} \\
  $5$ & 1.0 & 5 & 6136 & \num{6.03604} & \num{17.4764} & \num{0.0828331}
  \end{tabular}
  \caption{Iteration numbers, total energy, tangent-point potential and isometry errors of the final iterates in Example~\ref{ex:compress} for different choices of the TP exponent $q$ and refinement levels $k$ defining $\hat{h}_k = 2^{-k}$ and exponent $\b$ in $\rho = (\hat{h}/2)^\b$.}
  \label{tab:compress-table}
\end{table}

\begin{example}\label{ex:compress}
Let 
\[
\o=(-5,5)\times(0,1),\quad \GD = \{-5\} \times[0,1] \cup \{5\} \times[0,1],
\]
with boundary conditions 
\[
y_D(x_1,x_2) = [\a x_1, x_2, 0]^\top, \quad \phi_D = [I_2,0]^\top
\]
on $\GD$ for $\a = 0.1$, and let $f(x) = [0,0,c_f]^\top$ with $c_f = 10^{-6}$ for $x \in \o$.
We compute the resulting numerical equilibrium states on triangulations $\Th = \mathcal{T}_k$ consisting of halved squares with side lengths $\hat{h} = 2^{-k}$, $k=2,3,4,5$ with step size $\tau = \hat{h}/10$ for several values of the exponent $q$ and the TP parameter $\rho$ which is chosen as $\rho = (\hat{h}/2)^\b$, $\b = 0.5, 1, 1.5$.
Resulting numerical equilibrium configurations are depicted in Figure~\ref{fig:compress}.
Self-intersections are successfully prevented and the \emph{distance from self-contact} is decreasing for smaller values of $\rho$ corresponding to larger values of $\b$ and/or finer grids.
We observe that for higher exponents $q$ the tangent-point potential is more locally concentrated in potential contact regions whereas lower exponents lead to stronger repulsive effects throughout the whole domain.
The iteration numbers as well as the energies and isometry errors of the final iterates are listed in Table~\ref{tab:compress-table}. For the piecewise linear interpolants of the discrete solutions in the case $q = 5$, $\b = 1.0$, we obtain an experimental order of convergence of approximately~$\mathcal{O}(h^{1.6})$ in $L^2$.
\end{example}

\subsection{Compression of a twisted strip}
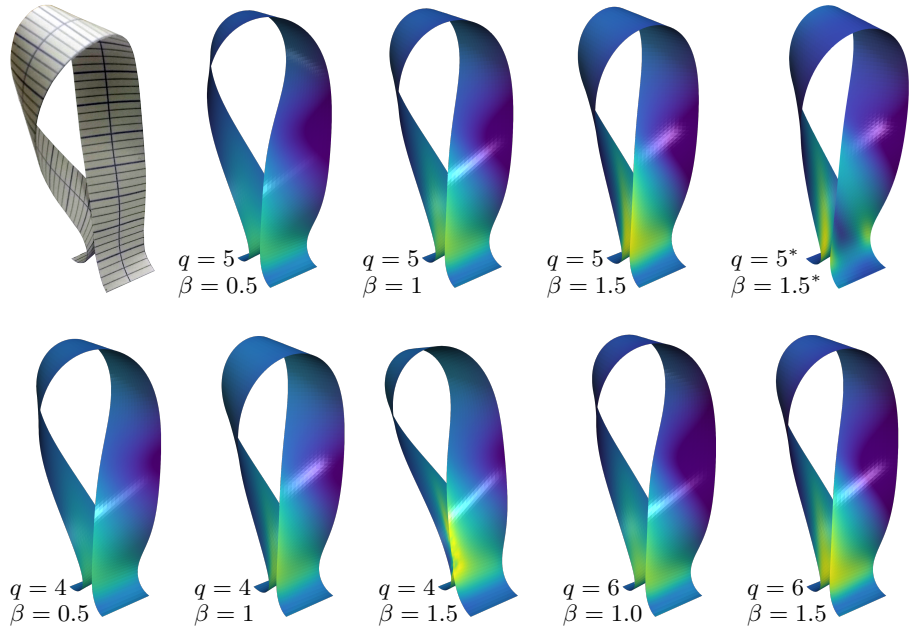
\begin{figure}
  \subfloat{\begin{overpic}[trim = -40mm -10mm -40mm -15mm, clip, height=40mm]{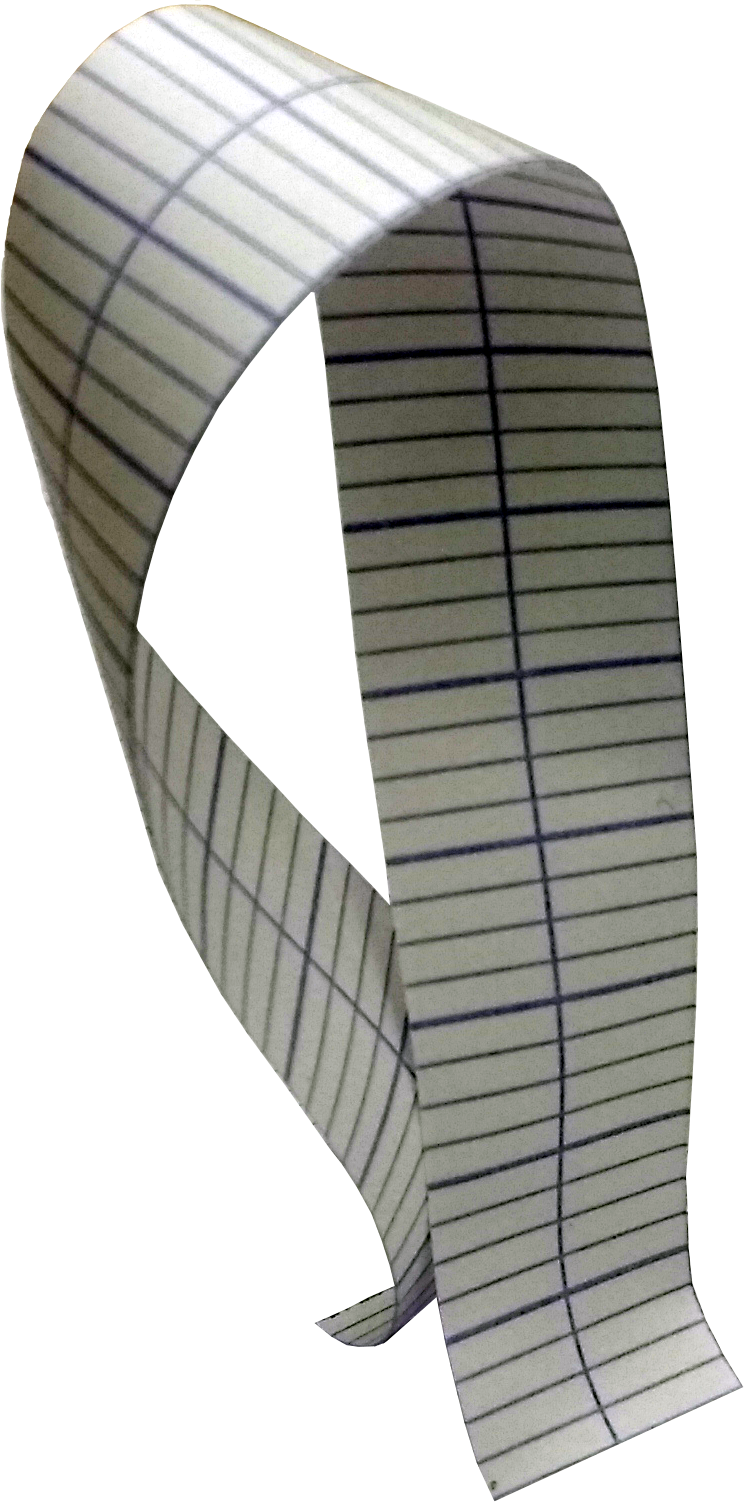}
  \end{overpic}}
  \subfloat{\begin{overpic}[trim = 60mm 5mm 20mm 0mm, clip, height=40mm]{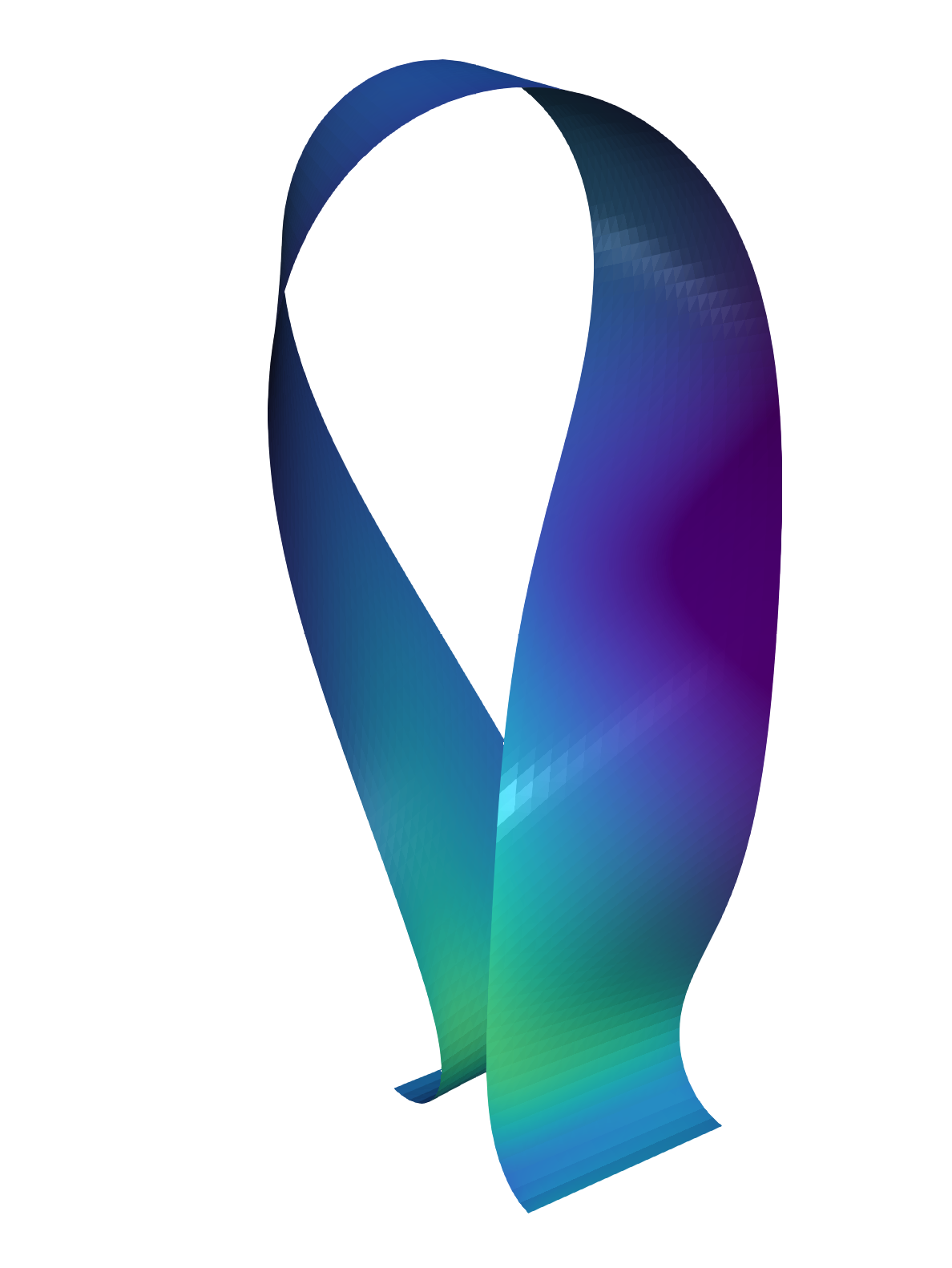}
  \put(1,2){\small $\b=0.5$}
  \put(1,12){\small $q=5$}
  \end{overpic}}
  \subfloat{\begin{overpic}[trim = 60mm 5mm 20mm 0mm, clip, height=40mm]{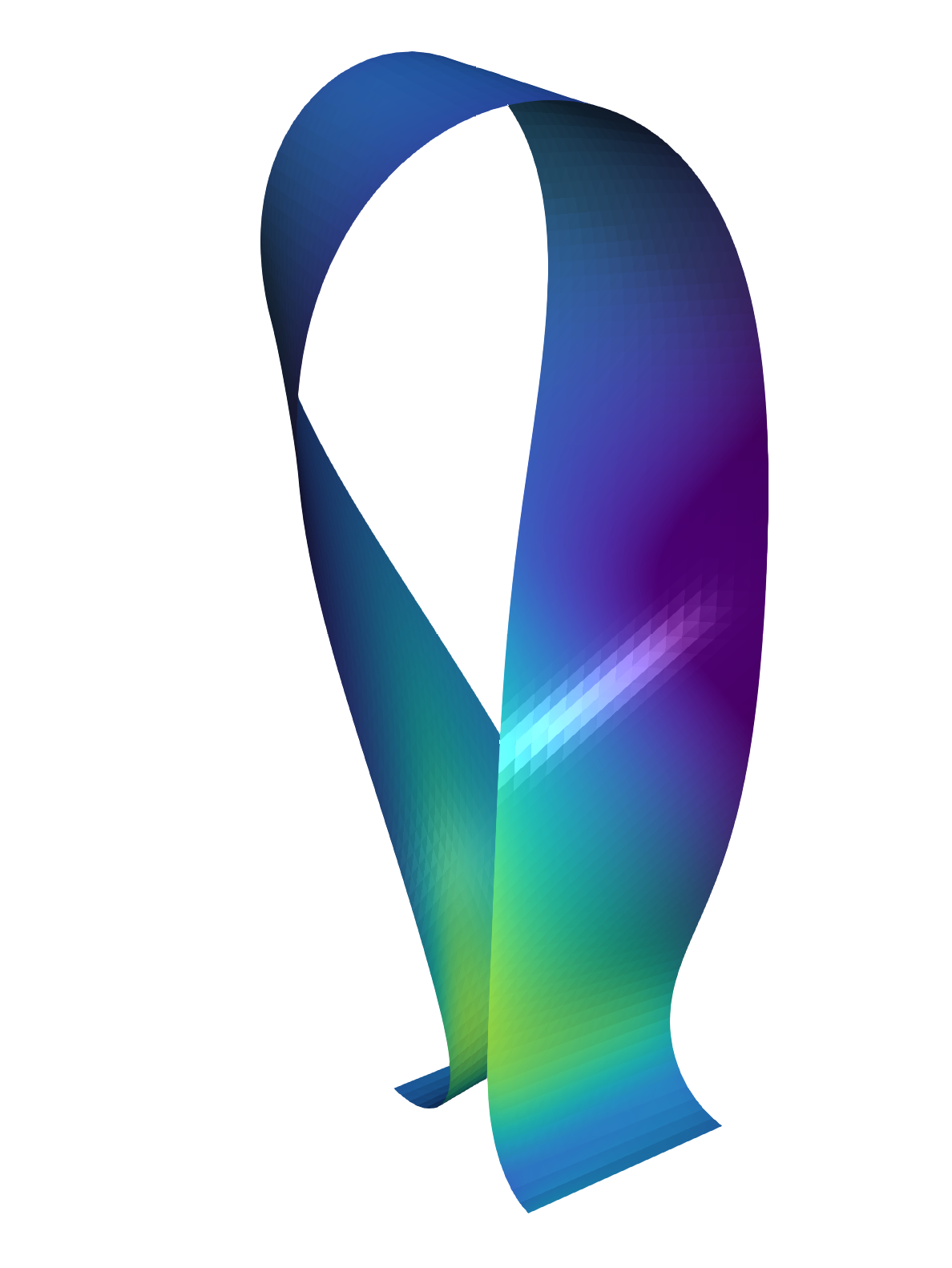}
  \put(1,2){\small $\b=1$}
  \put(1,12){\small $q=5$}
  \end{overpic}}
  \subfloat{\begin{overpic}[trim = 60mm 5mm 20mm 0mm, clip, height=40mm]{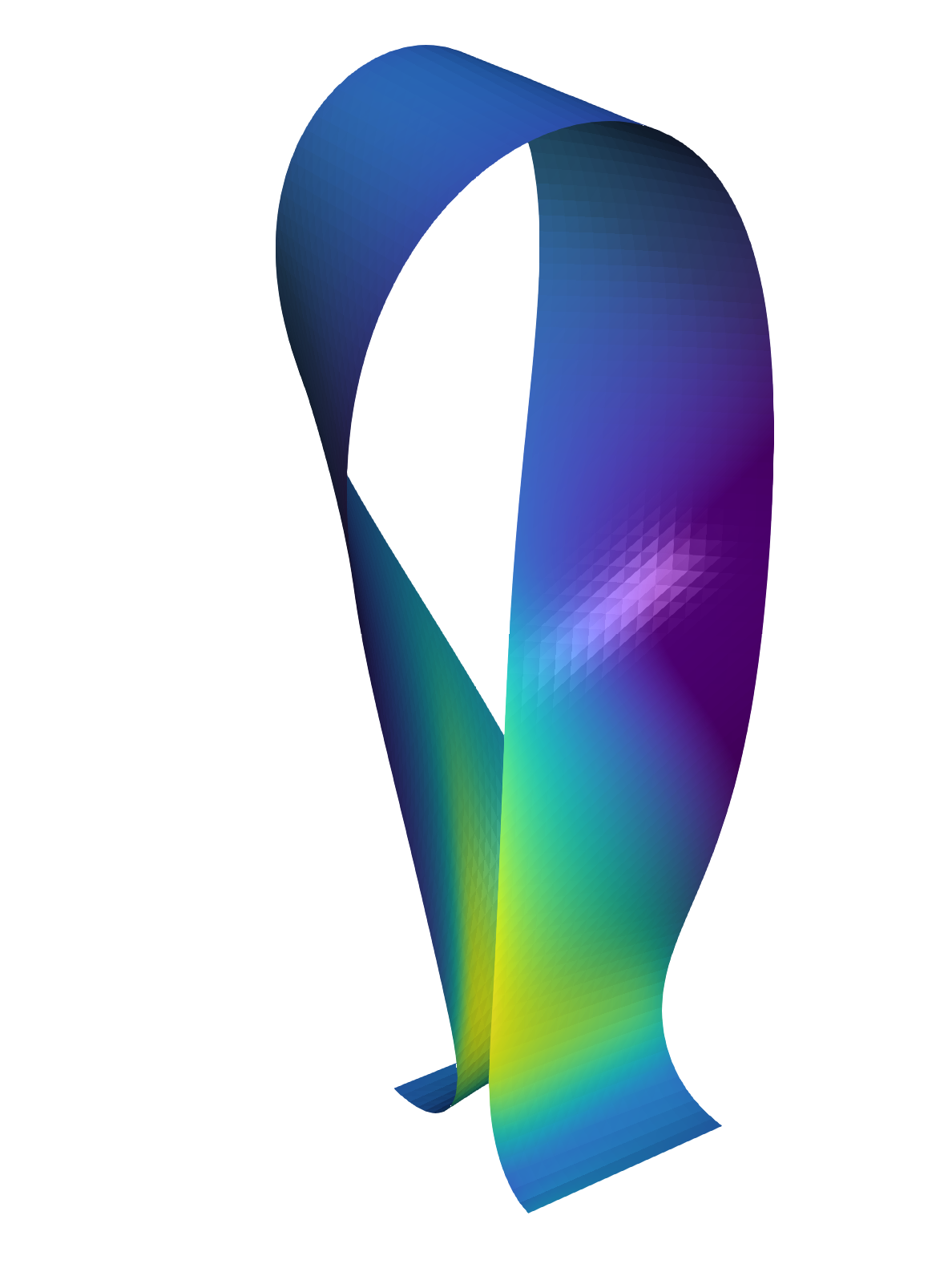}
  \put(1,2){\small $\b=1.5$}
  \put(1,12){\small $q=5$}
  \end{overpic}}
  \subfloat{\begin{overpic}[trim = 60mm 5mm 20mm 0mm, clip, height=40mm]{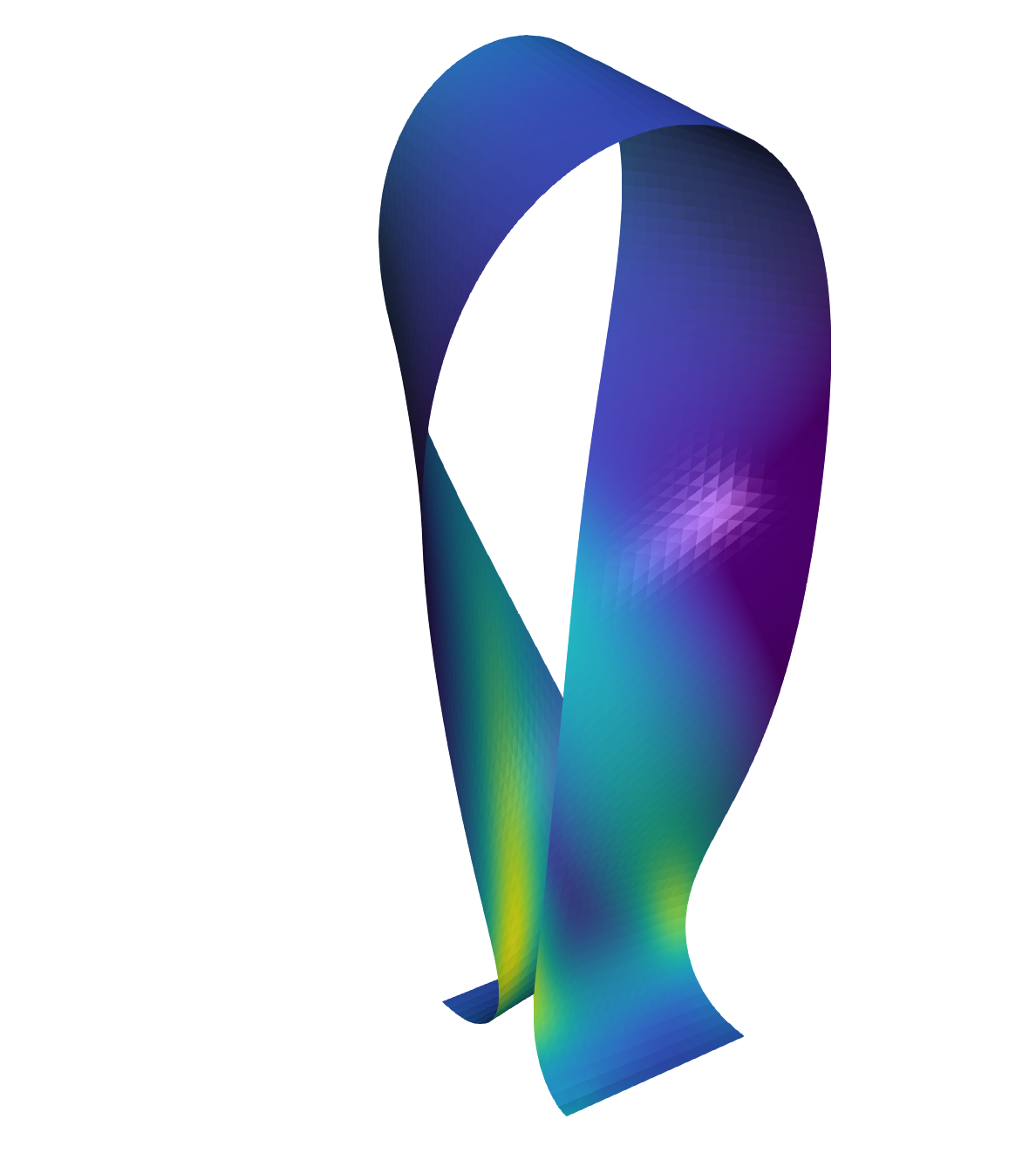}
  \put(1,2){\small $\b=1.5^*$}
  \put(1,12){\small $q=5^*$}
  \end{overpic}}\\
  % newline
  \subfloat{\begin{overpic}[trim = 60mm 5mm 20mm 0mm, clip, height=40mm]{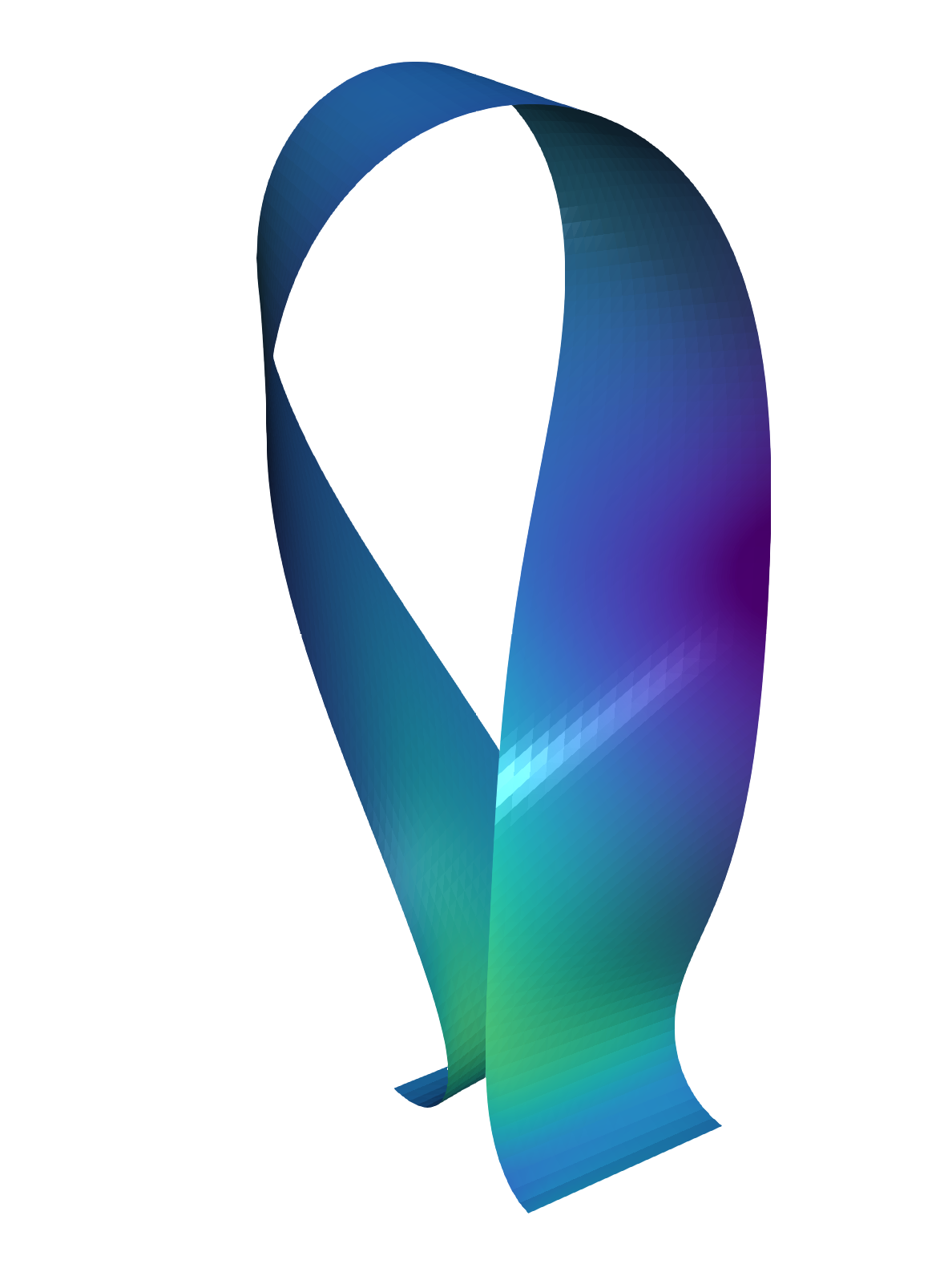}
  \put(1,2){\small $\b=0.5$}
  \put(1,12){\small $q=4$}
  \end{overpic}}
  \subfloat{\begin{overpic}[trim = 60mm 5mm 20mm 0mm, clip, height=40mm]{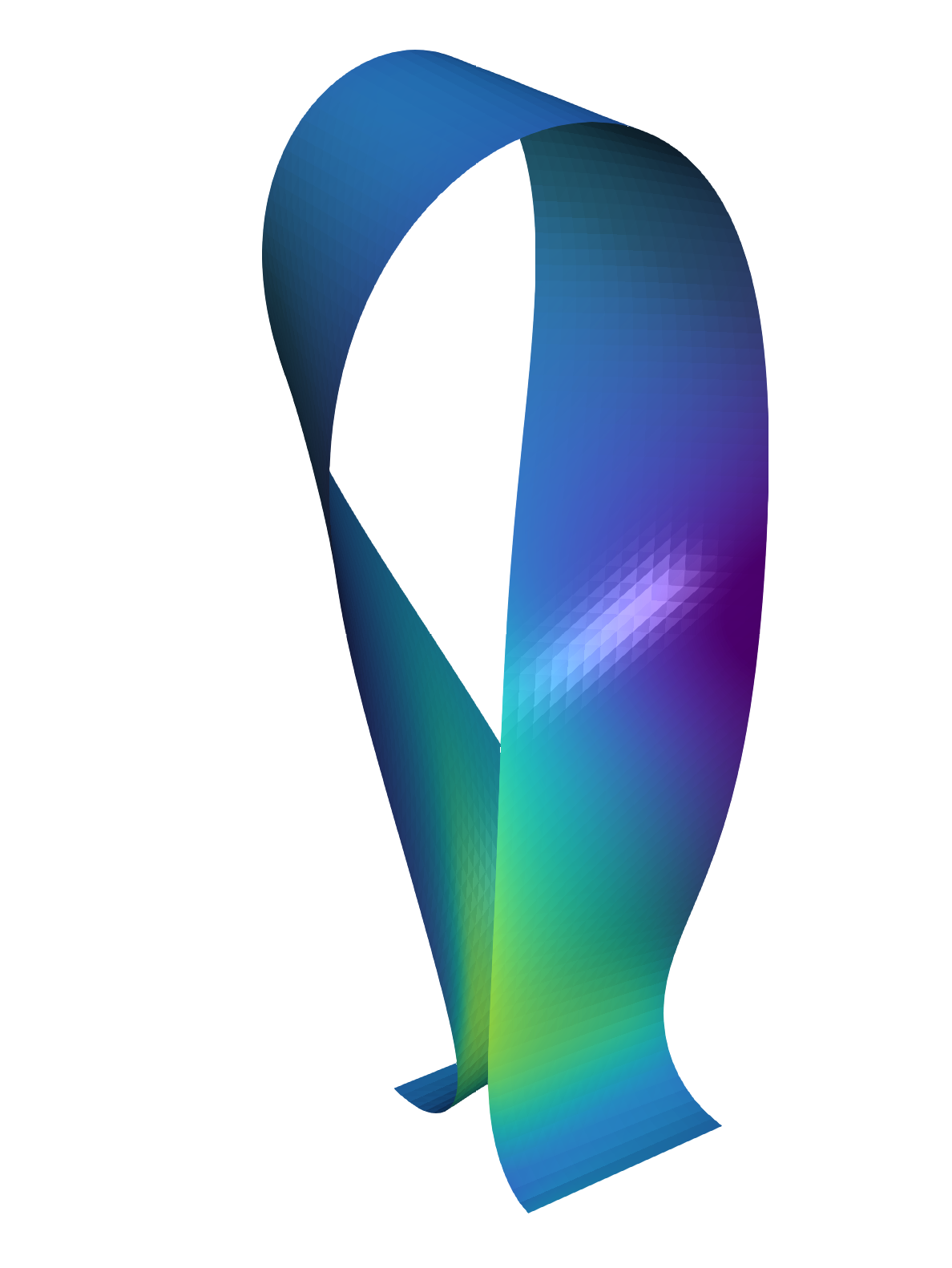}
  \put(1,2){\small $\b=1$}
  \put(1,12){\small $q=4$}
  \end{overpic}}
  \subfloat{\begin{overpic}[trim = 60mm 5mm 20mm 0mm, clip, height=40mm]{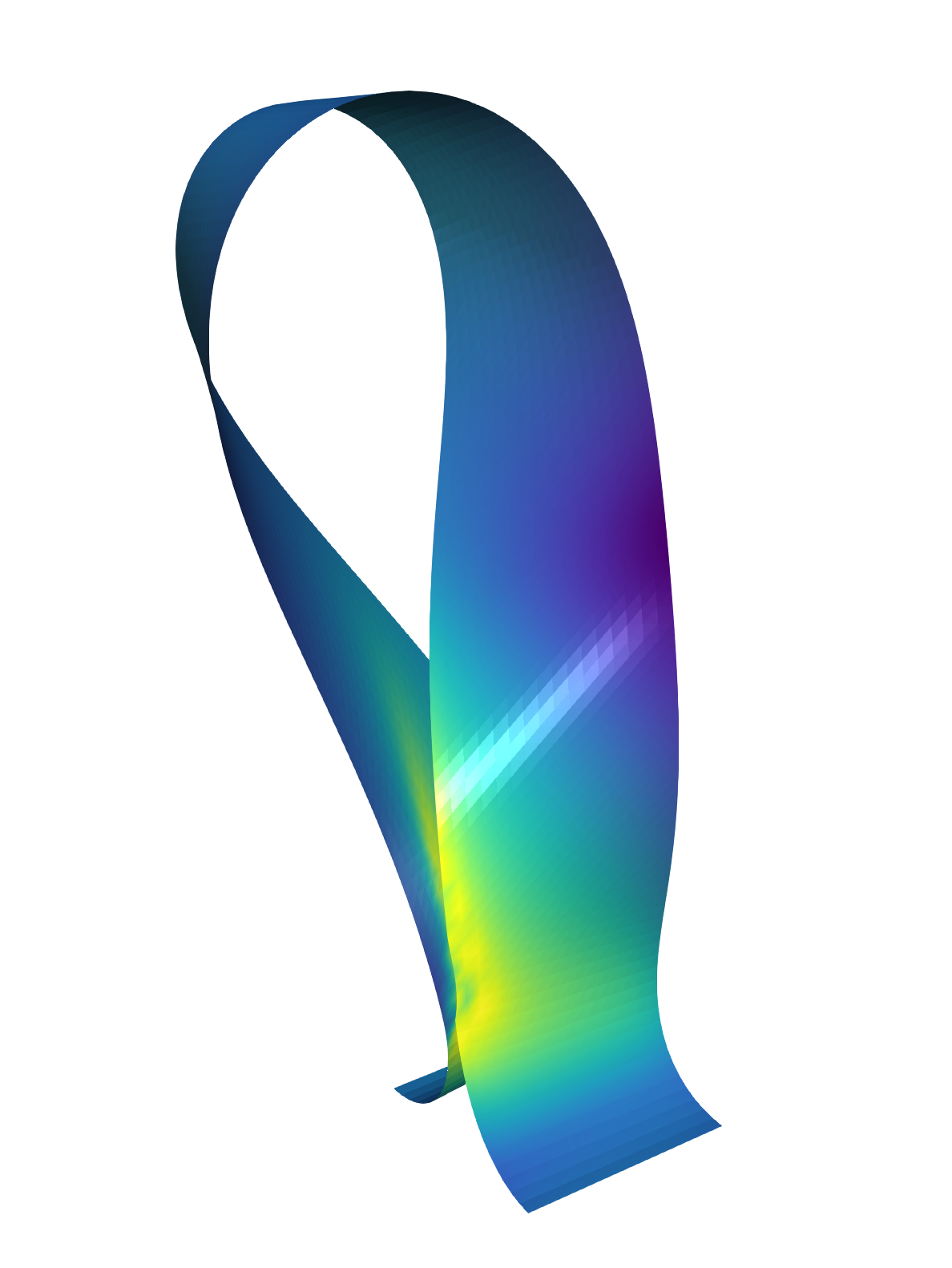}
  \put(1,2){\small $\b=1.5$}
  \put(1,12){\small $q=4$}
  \end{overpic}}
  \subfloat{\begin{overpic}[trim = 60mm 5mm 20mm 0mm, clip, height=40mm]{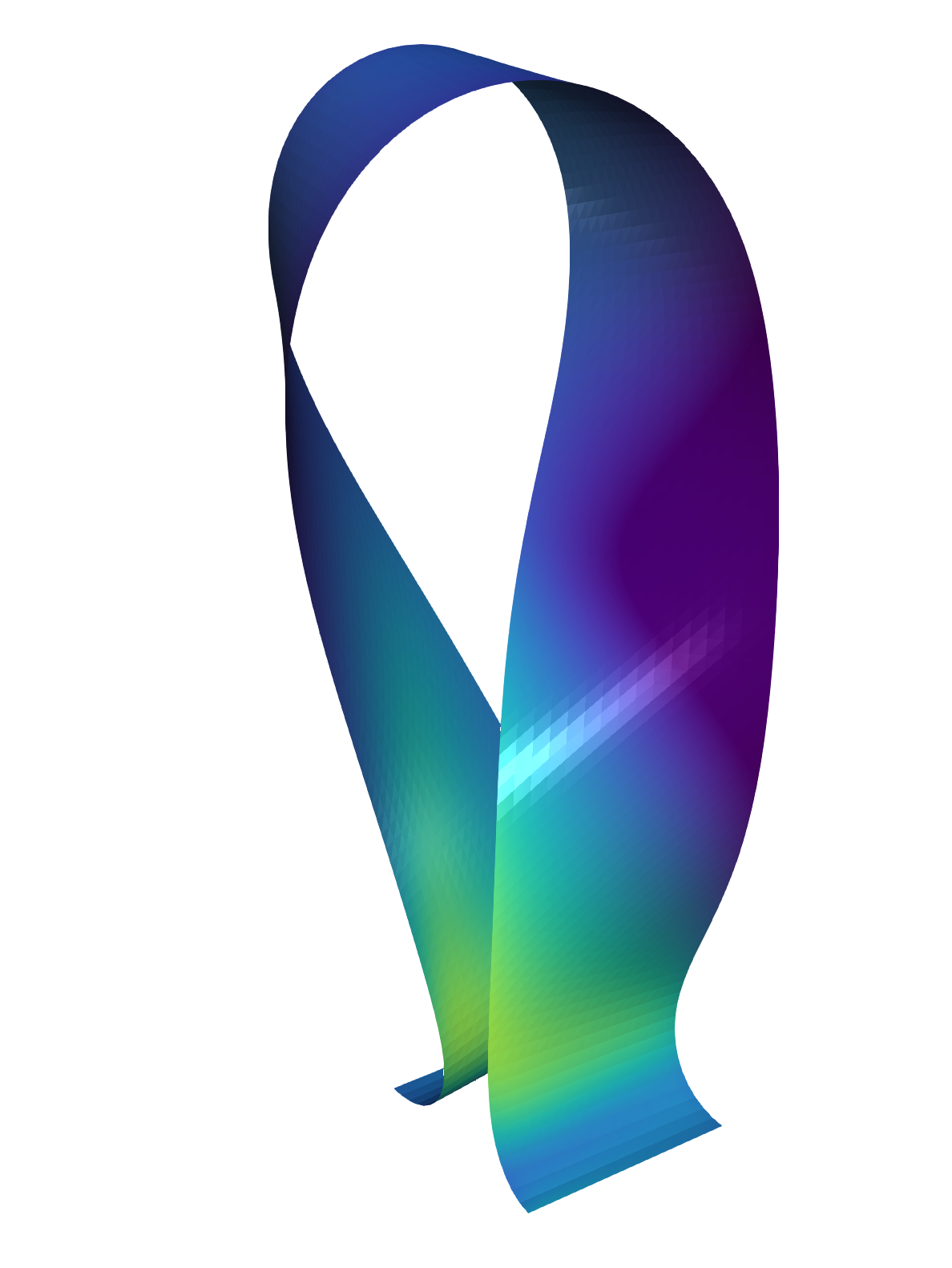}
  \put(1,2){\small $\b=1.0$}
  \put(1,12){\small $q=6$}
  \end{overpic}}
  \subfloat{\begin{overpic}[trim = 60mm 5mm 20mm 0mm, clip, height=40mm]{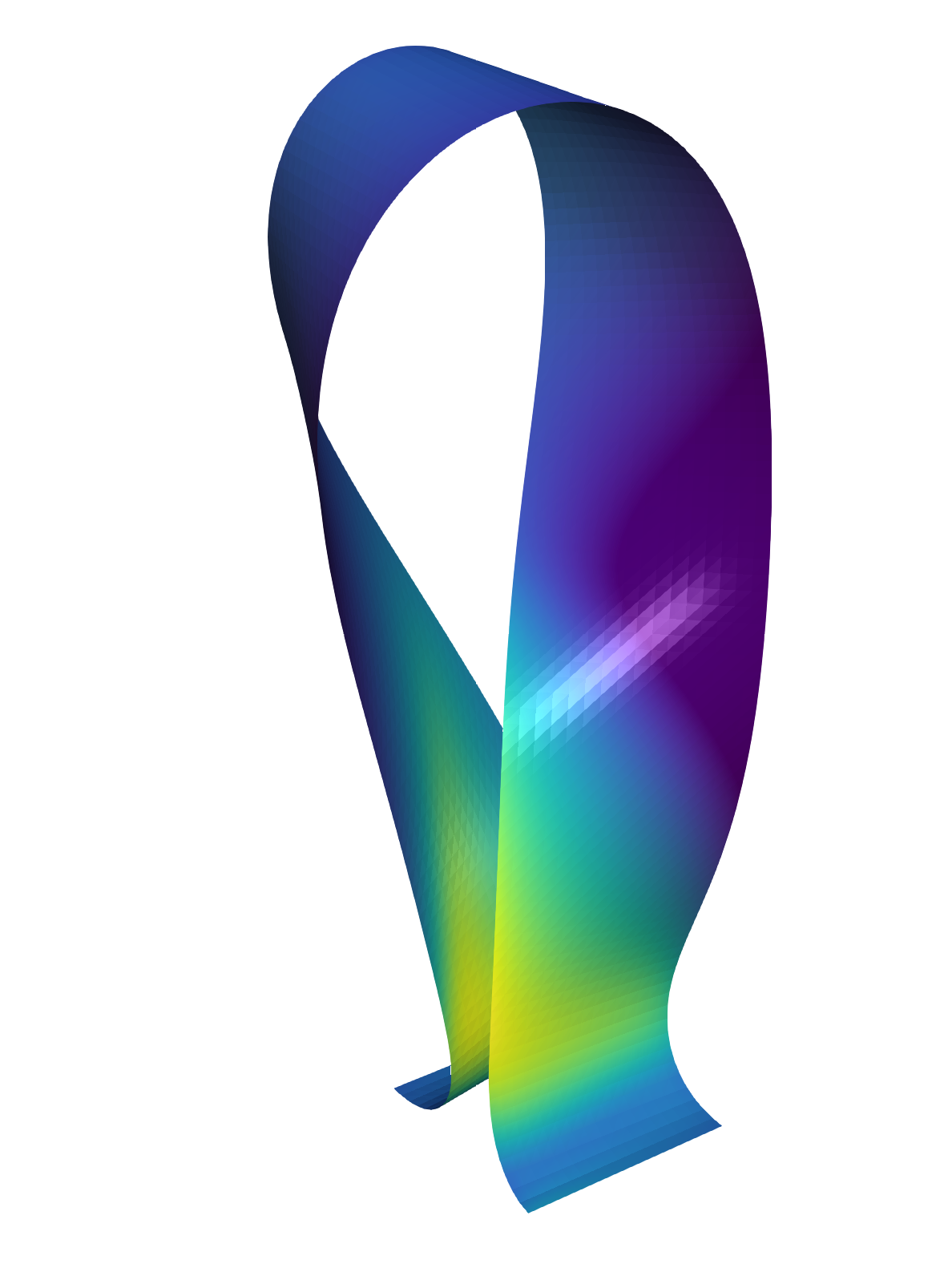}
  \put(1,2){\small $\b=1.5$}
  \put(1,12){\small $q=6$}
  \end{overpic}}  \\
  \caption{Photograph of an actual paper strip (top left) and stationary configurations of the compressed twisted strip in Example~\ref{ex:twist180} for $\#\Th = 5120$ and different values of the exponent $\b$ defining the self-avoidance parameter $\rho = (\hat{h}/2)^\b$ and TP exponent~$q$.
  Note how smaller values of~$q$ imply a stronger repulsive effect throughout the whole domain. As expected the numerical solutions corresponding to smaller values of $\rho$ and higher values of $q$ bear a closer resemblance to the photograph. Self-intersection is not prevented in the case $q=4, \beta=1.5$. The plot marked with an asterisk on the top right corresponds to the stationary configuration resulting from a replacement of the tangent-point potential with a discrete boundary-domain potential.}
  \label{fig:twist180}
\end{figure}

\begin{figure}
  \centering
  \includegraphics[width=0.7\linewidth]{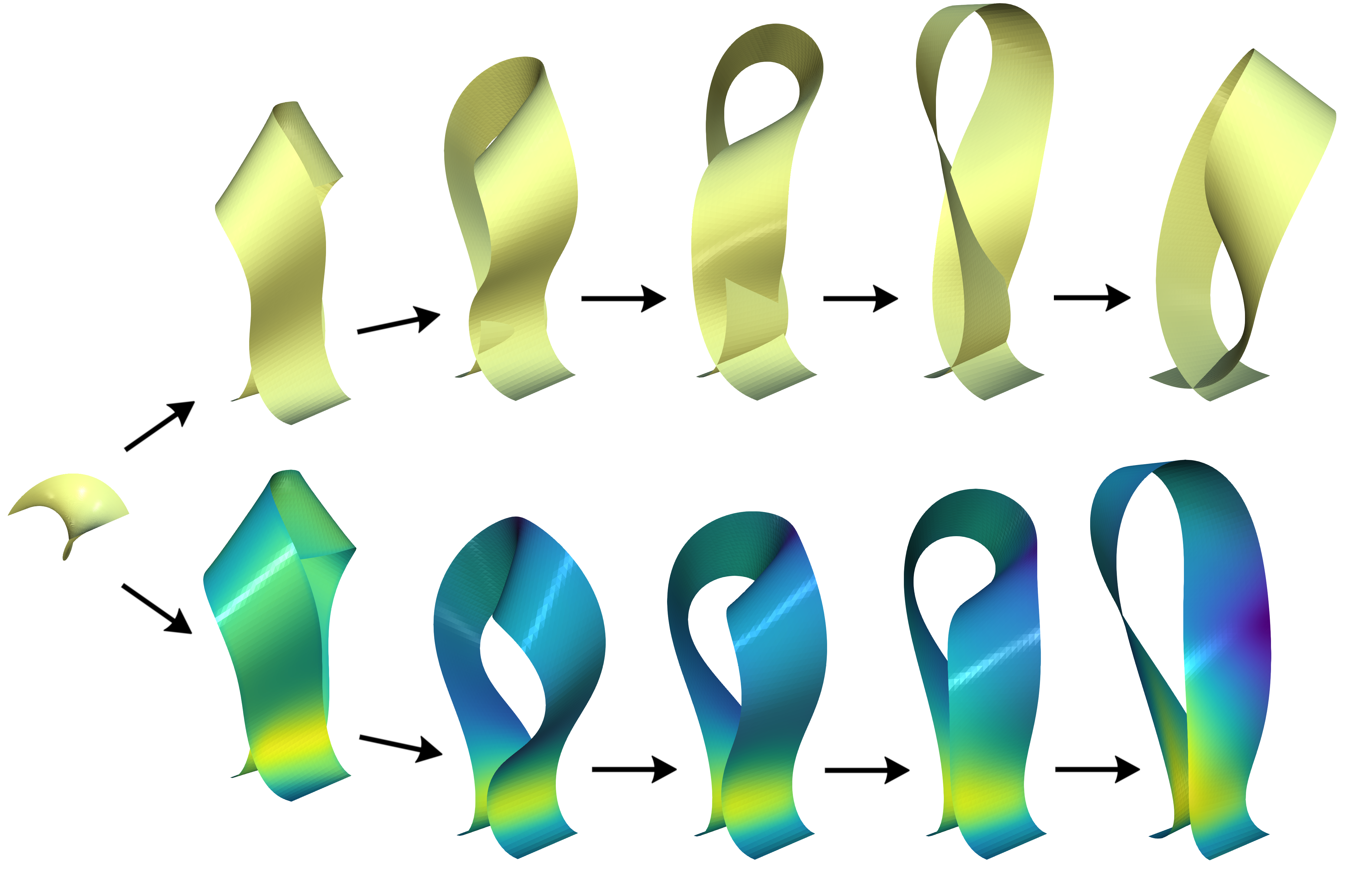}
  \put(-290,120){$\rho = 0$}
  \put(-290,32){$\rho > 0$}
  \caption{Comparison of the evolution of a self-avoiding strip versus a non-self-avoiding strip subject to compressive boundary conditions implying a twist of 180 degrees in Example~\ref{ex:twist180}.}
  \label{fig:twist180evol}
\end{figure}

For the same reference configuration of a $10 \times 1$ plate as in Example~\ref{ex:compress} we modify the boundary conditions and initial value to model a 180-degree twist of the strip in addition to its compression. 
As in Example~\ref{ex:compress}, the numerically computed minimizers of the pure bending energy (corresponding to the choice $\rho=0$) exhibit self-intersections, cf. Figure~\ref{fig:twist180evol} where we compare the evolution of the discrete gradient flow in the cases $\rho=0$ and $\rho>0$. 
Furthermore, we presume that self-contact in corresponding physical solutions is more singular than the expected self-contact in Example~\ref{ex:compress} in the sense that the (almost-)contact region now may contain isolated one-dimensional subsets of the domain boundary, cf. Figure~\ref{fig:twist180}.
\begin{table}
  \centering
  \begin{tabular}{c c c c c c c}\hline
  $k$ & $\b$ & $q$ & $N_\mathrm{iter}$ & $E_h[y_h^\infty]$ & $\TP_{h}[y_h^\infty]$ & $\delta_\mathrm{iso}[y_h^\infty]$ \\ \hline
  $2$ & 1.0 & 5 & 1315  & \num{8.93104} & \num{5.55138}  & \num{0.228221} \\
  $3$ & 1.0 & 5 & 4367  & \num{8.94381} & \num{8.93415}  & \num{0.0965555} \\
  $4$ & 0.5 & 4 & 11718 & \num{9.64297} & \num{5.38592}  & \num{0.0726214} \\
  $4$ & 1.0 & 4 & 12454 & \num{8.4946} &  \num{14.0203}  & \num{0.07224} \\
  $4$ & 1.5 & 4 & 10196 & \num{8.19601} & \num{43.1725}  & \num{0.0722287} \\
  $4$ & 0.5 & 5 & 11469 & \num{10.1082} & \num{4.84628}  & \num{0.113205} \\
  $4$ & 1.0 & 5 & 12220 & \num{8.91663} & \num{14.4271}  & \num{0.0725601} \\
  $4$ & 1.5 & 5 & 12683 & \num{8.31079} & \num{46.5317}  & \num{0.0721425} \\
  $4$ & 0.5 & 6 & 51012 & \num{22.7848} & \num{8.17345}  & \num{2.4617} \\
  $4$ & 1.0 & 6 & 12124 & \num{9.51322} & \num{14.9389}  & \num{0.169585} \\
  $4$ & 1.5 & 6 & 12648 & \num{8.6628}  & \num{51.5257}  & \num{0.0774521} \\
  $4^*$ & 1.0$^*$ & 5$^*$ & 13562$^*$ & \num{8.50455}$^*$ & \num{54.3076}$^*$ & \num{0.0722385}$^*$ \\
  \end{tabular}
  \caption{Iteration numbers, total energy, tangent-point potential and isometry error of the final iterates in Example~\ref{ex:twist180}. Note, how in the case $k=4$, $\b= 0.5$, $q=6$, the repulsive force is locally too strong for the chosen step size, resulting in a more severe violation of the discrete isometry constraint. The values in the last line were obtained using a discrete version of the boundary-domain potential $\TP_{\p \o, \veps}$ instead of $\TP_{h}$. Consequently, the value in the $\TP_{h}$ column in the last line corresponds to the discretized boundary-domain potential.  }
  \label{tab:twist180-table}
\end{table}
\begin{figure}
  \centering
  \includegraphics[width=.9\linewidth]{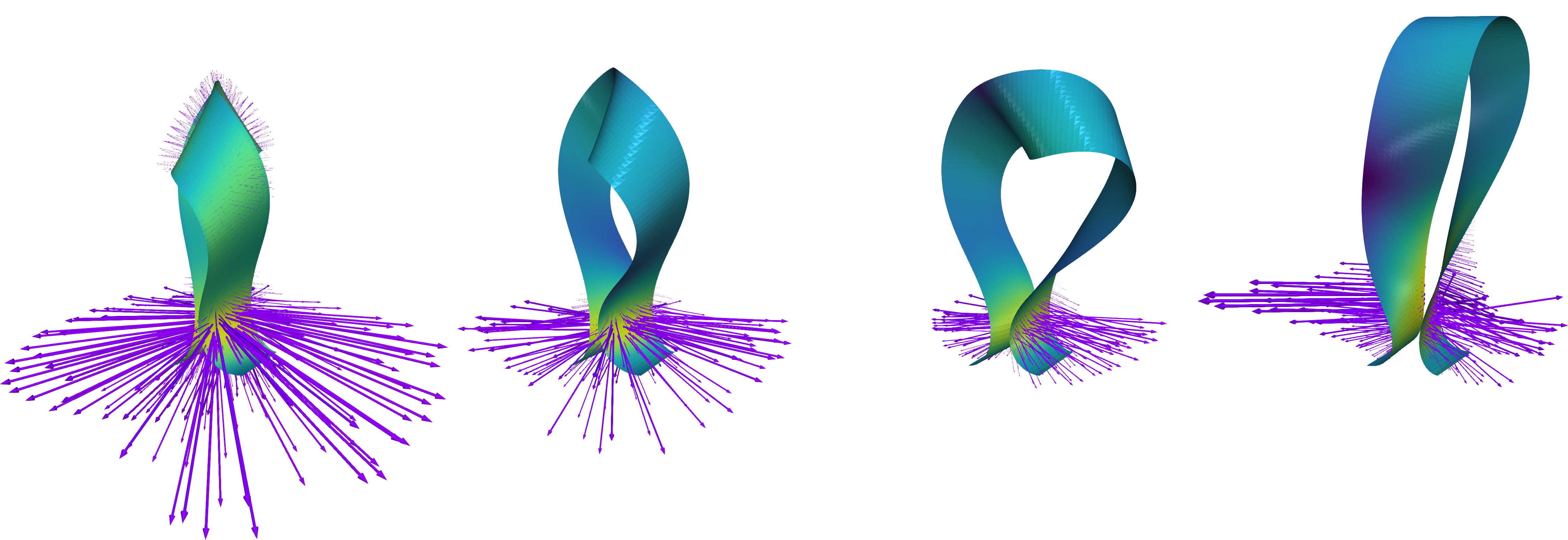}
  \caption{Visualization of the repulsive pseudo-force field $f_\tp$ which is induced by the tangent-point potential in several configurations of the twisted $10 \times 1$ strip with compressive boundary conditions in Example~\ref{ex:twist180}.}
  \label{fig:tpforce}
\end{figure}
\begin{figure}
  \centering
  \includegraphics[width=.84\textwidth]{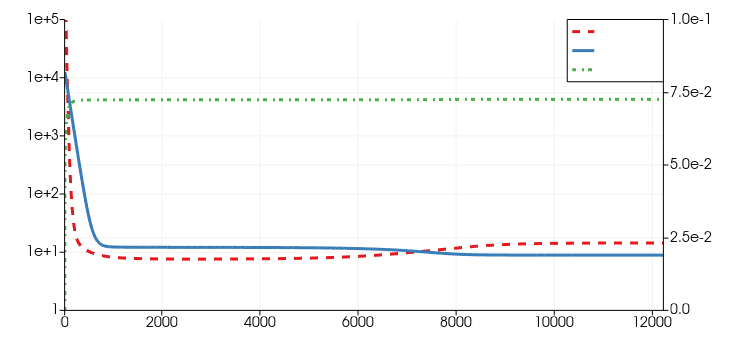}
  \put(-182,1){\small $k$}
  \put(-372,95){\small $E_h[y_h^k],$}
  \put(-382,83){\small $\mathrm{TP}_{\veps,h}[y_h^k]$}
  \put(-4.5,88){\small $\delta_\mathrm{iso}[y_h^k]$}
  \put(-61,151.5){\footnotesize $\mathrm{TP}_{\veps,h}$}
  \put(-61,142){\footnotesize $E_h$}
  \put(-61,133){\footnotesize $\delta_\mathrm{iso}$}
  \caption{Evolution of the tangent-point potential, total energy decay and boundedness of the isometry error for the case $q = 5$, $\beta=1$ in Example~\ref{ex:twist180}. The decrease of $E_h[y_h^k]$ between $k=6000$ and $k=8000$ corresponds to the qualitative change of the almost-contact from ``surface-to-surface'' to ``edge-to-surface'', cf. the second last versus the last snapshot of the evolution for $\rho>0$ depicted in Figure~\ref{fig:twist180evol} as well as the second last and last configuration in Figure~\ref{fig:tpforce}.}
  \label{fig:twist180energy}
\end{figure}

\begin{example}\label{ex:twist180}
  We let $\o = (-5,5)\times(0,1)$, $\GD = \{-5\} \times[0,1] \cup \{5\} \times[0,1] $, $f(x) = [0,0,c_f]^\top$ with $c_f = 10^{-6}$ for $x \in \o$, 
  \[
  y_D(x_1,x_2) = \begin{cases}
    [\a x_1, x_2, 0]^\top, &\text{ if } x_1 = -5, \\
    [\a x_1, 1 - x_2, 0]^\top, &\text{ if } x_1 = 5,
  \end{cases}
  \]
  with $\a = 0.1$, and $\phi_D = \left[\operatorname{diag}\left(1,-\sign(x_1)\right),0\right]^\top$ for $(x_1,x_2) \in \GD$.
  For the iteration we choose a simple extension of the boundary data as initial value which does not correspond to an actual isometry in $H^2$. 
  In order to overcome the difficulties arising from steep gradients of the tangent-point potential of the initial value, we set $\rho$ to zero and perform some iterations to relax the initial data. 
  We then set $\rho$ to a positive value and restart the discrete gradient flow using the relaxed data as the initial value. 
  Figure~\ref{fig:twist180} shows a photograph of an actual paper strip in comparison with the calculated numerical equilibrium states for a triangulation consisting of halved squares of side length $\hat{h} = 2^{-4}$ amounting to 5120 triangles. 
  The step size was chosen as $\tau = \hat{h}/10$ and the final configurations correspond to different values of the exponent $q$ and the TP parameter $\rho = (\hat{h}/2)^\b$. In the case $q=6$, $\beta=0.5$ the occurrence of strong repulsive pseudo forces in our heuristically computed initial value immediately lead to a comparatively large violation of the isometry constraint, cf. Table~\ref{tab:twist180-table}. For smaller values of $q$ stronger repulsive effects of the tangent-point potential throughout the whole domain can be observed. In the case $q=4$, $\b=1.5$ the repulsive effects are not strong enough to successfully prevent self-intersection and the resulting deformation is not injective. 
  As one would expect, numerical solutions corresponding to smaller values of $\rho$ and higher values of $q$ bear a closer resemblance to the photograph, at least in the eyeball metric. For all investigated choices of parameters $\b, q$, we observed a tendency towards the physical configuration for $h \to 0$, which was more pronounced for higher values of the TP exponent $q$. This is, however, a purely qualitative comparison, as no material parameters have been determined.
  Corresponding iteration numbers, energies, tangent-point potentials and isometry errors of the final iterates are listed in Table~\ref{tab:twist180-table}.
  The last line of this table, marked with an asterisk, contains values that are obtained when the tangent-point potential $\TP_h$ is replaced with a discretization of the boundary-domain potential $\TP_{\p \o, \veps}$ while all other parameters remain unchanged. The slightly increased number of necessary iterations is more than compensated for by the faster numerical integration of this potential. For comparison a plot of the corresponding final configuration is shown in Figure~\ref{fig:twist180}, also marked with an asterisk there. The stronger concentration of the boundary-domain potential around regions of possible contact points seems to indicate that its use is more justifiable from a physical point of view in this experiment.
  
  For the case $q=5$ and $\b = 1$ we compare the evolution of the discrete gradient flow with and without tangent-point potential in Figure~\ref{fig:twist180evol}. 
  For the same case the induced pseudo-force field $f_\tp$ is visualized for several configurations in Figure~\ref{fig:tpforce}. 
  The boundedness of the isometry error and the monotone decay of the total energy despite an increase in the tangent-point potential are illustrated in Figure~\ref{fig:twist180energy}. For the piecewise linear interpolants of the discrete solutions, we obtain an experimental order of convergence in $L^2$ of approximately~$\mathcal{O}(h^{0.9})$. However, due to the higher iteration numbers in this experiment the employed meshes are coarser than in the previous example and the asymptotic range may not have been reached, yet.
\end{example}

\subsection{Effects of torsion in a circular ribbon}
We investigate the effects of imposing torsion on the energy minimizing configurations of a periodic ribbon. To this end we prescribe closed loop boundary conditions on the short ends of a $50 \times 1$ ribbon and choose the initial value for the discrete gradient flow such that the ribbon contains a number $K$ of 180-degree twists, e.\,g. $K=0$ corresponds to a simple closed loop and $K=1$ to a Möbius strip.
\begin{figure}
  \centering
  \includegraphics[width=0.94\textwidth]{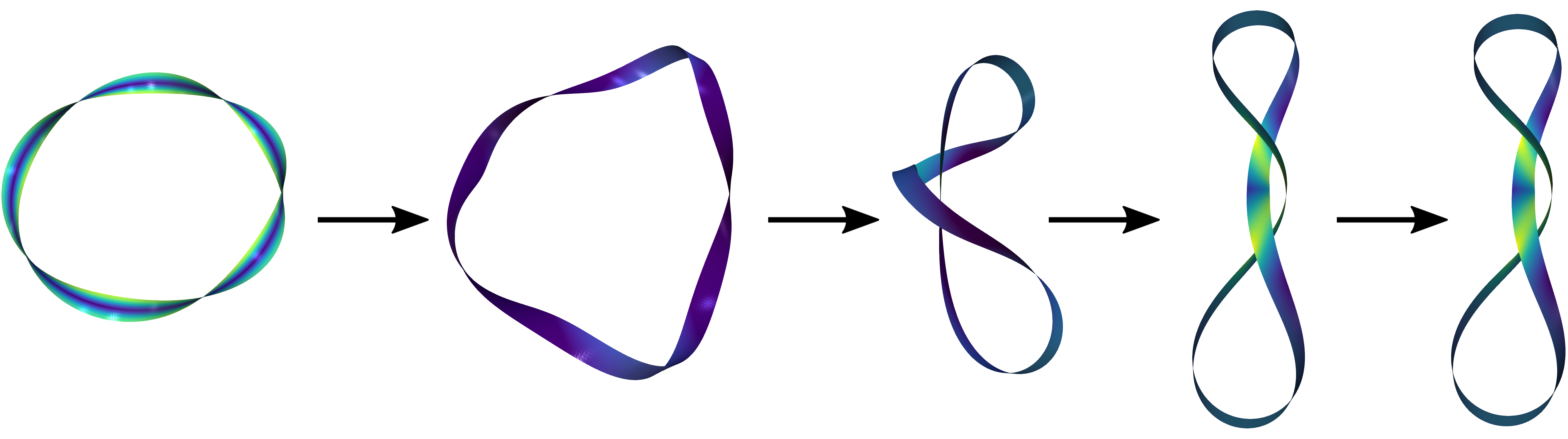}
  \caption{Snapshots of the evolution for the closed twisted ribbon in Example~\ref{ex:twist_circ} after 0, 300, 500, 750, 2000 and 3463 iterations.}
  \label{fig:twisted_circle_evo}
\end{figure}

\begin{example}\label{ex:twist_circ}
  Let $\o=(0,50)\times(0,1)$, $f = [0,0,0]^\top$ and $K=5$. We employ a triangulation of $\o$ into 6400 triangles given by halved squares with side length $\hat{h} = 2^{-3}$ and define the preliminary initial value $\ty_h^0$ via
  \[ 
  \ty_h^0(z) = \begin{bmatrix}
    6 \cos(2\pi z_1/50) \\ 6 \sin(2\pi z_1/50) \\ 0
  \end{bmatrix} + \sin(K\pi z_1/50) \begin{bmatrix}
    \cos(2\pi z_1/50) \\ \sin(2\pi z_1/50) \\ 0
  \end{bmatrix} + (z_2-0.5) \begin{bmatrix}
    0 \\ 0 \\ \cos(K\pi z_1/50)
  \end{bmatrix}, 
  \]
  and
  \[ 
  \nabla \ty^0_h(z) = \begin{bmatrix}
    -\sin(2\pi z_1/50) & \cos(2\pi z_1/50)\sin(K\pi z_1/50)\\
    \cos(2\pi z_1/50) & \sin(2\pi z_1/50)\sin(K\pi z_1/50)\\
    0 & \cos(K\pi z_1/50)
  \end{bmatrix}, 
  \]
  for all vertices $z = (z_1,z_2)\in\cN_h$. As in the previous example this initial data does not resemble an isometry in $H^2$. 
  However, the discrete isometry constraint is satisfied and, as before, we obtain the actual initial data $y_h^0$ for the algorithm from a relaxation of $\ty_h^0$.
  With step size $\tau = \hat{h}/10$, TP parameter $\rho = \hat{h}$ and exponent $q = 5$ the iteration terminates after $3463$ steps at an intertwined configuration with final energy $E_h[y_h^\infty] = 2.589$ and isometry error $\d_\mathrm{iso}[y_h^\infty] = 1.550 \times 10^{-2}$. Snapshots of the corresponding evolution are depicted in Figure~\ref{fig:twisted_circle_evo}.
\end{example}

\subsection{Trefoil knot}\label{subsec:trefoil}
As in Example~\ref{ex:twist_circ} we consider a $50 \times 1$ strip with closed boundary conditions, but now choose the initial data such that it bears the topology of a trefoil knot, which we expect to be preserved by the algorithm.
\begin{figure}
  \subfloat{
    \includegraphics[height=2.8cm]{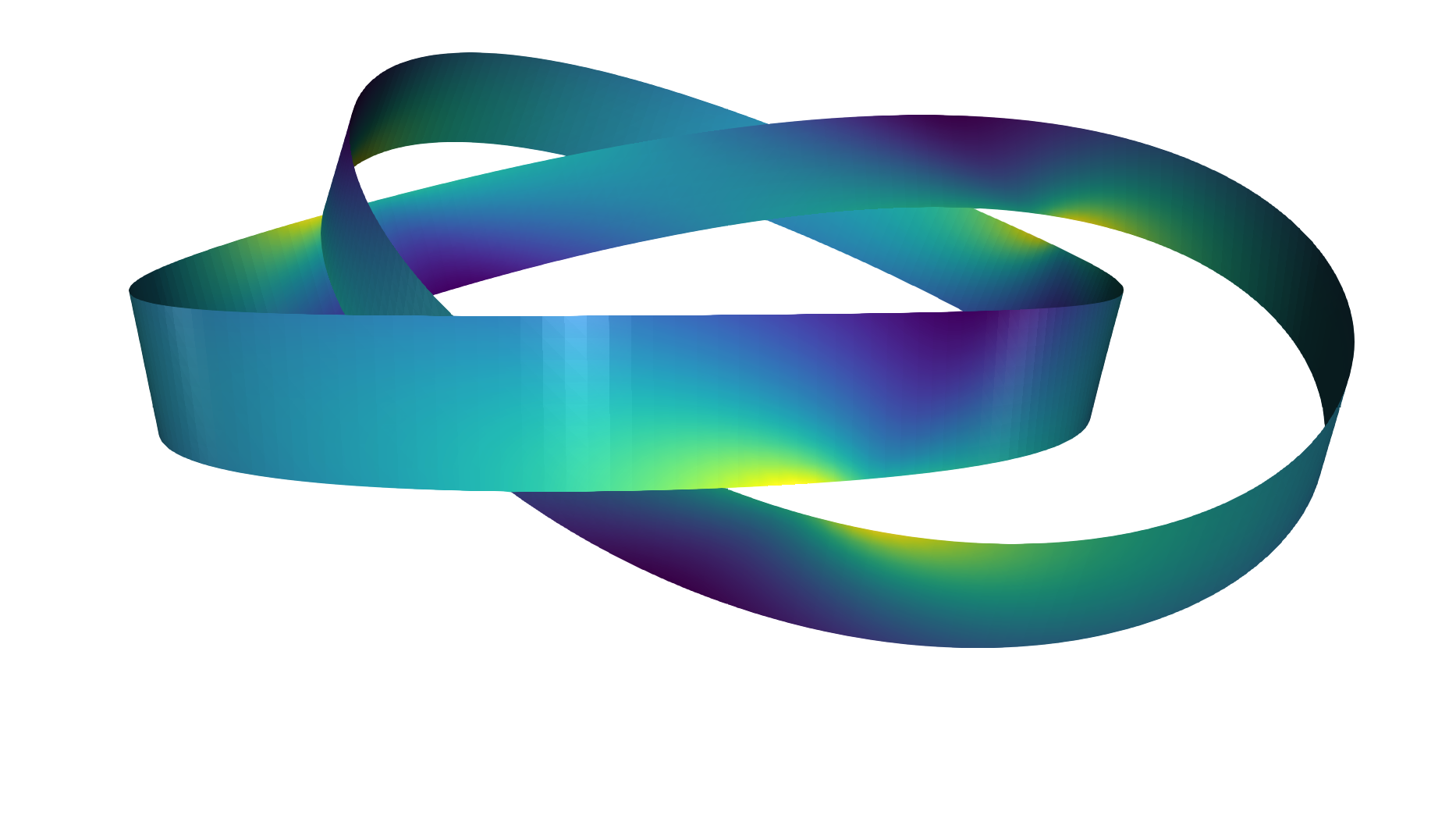}
  }
  \hspace{2mm}
  \subfloat{
    \includegraphics[height=2.8cm]{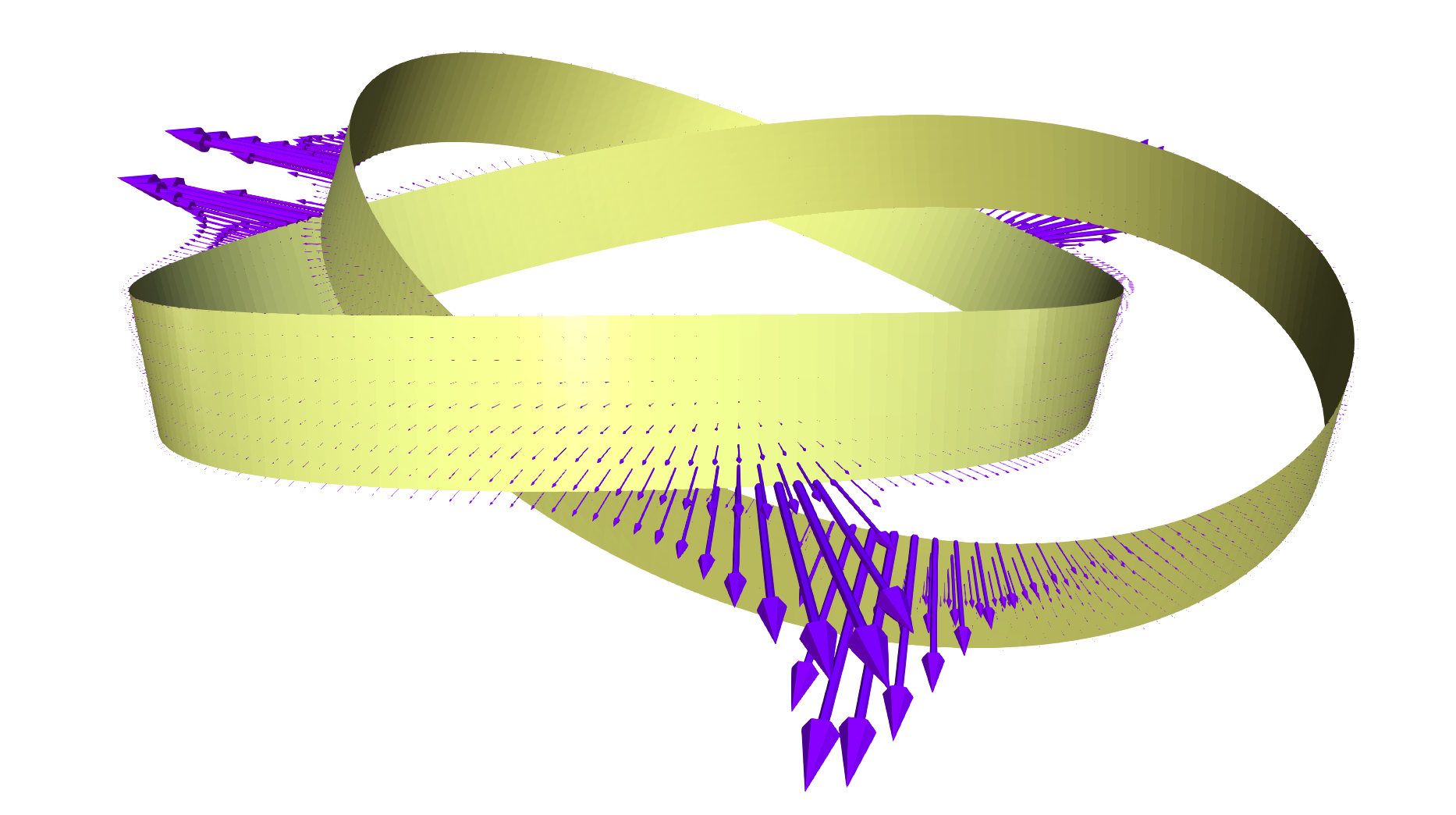}
  }
  \caption{
  Knot topology of the trefoil in Example~\ref{ex:trefoil} is preserved. Both plots show the final iterate in the case $k=3$. The glyphs in the right-hand-side plot are a visualization of the pseudo-force field $f_\tp$.
  }
  \label{fig:trefoil}
\end{figure}
\begin{table}
  \centering
  \begin{tabular}{c c c c c c c}\hline
    $k$ & $\tau$ &$\rho$ & $N_\mathrm{iter}$ & $E_h[y_h^\infty]$ & $\TP_{h}[y_h^\infty]$ & $\delta_\mathrm{iso}[y_h^\infty]$ \\ \hline
    $1$ & \num{0.0100}   & \num{0.25}   & 1078 & \num{2.04388} & \num{0.601009}   & \num{0.00197947} \\
    $2$ & \num{0.0050}   & \num{0.125}  & 2138 & \num{1.96454} & \num{0.911468} & \num{0.000866212} \\
    $3$ & \num{0.0025}  & \num{0.0625} & 4504 & \num{1.90316} & \num{1.37857}  & \num{0.00040601} \\
  \end{tabular}
  \caption{Iteration numbers, total energy, tangent-point potential and isometry error of the final iterates in the trefoil Example~\ref{ex:trefoil}.}
  \label{tab:trefoil-table}
\end{table}

\begin{example}\label{ex:trefoil}
  Let $\o=(0,50)\times(0,1)$ and $f = [0,0,0]^\top$. We consider the (one-dimensional) parametrization $u \colon [0,50] \to \R^3$ of a trefoil knot, 
  \begin{equation*}
    u(t) = \begin{bmatrix} u_1(t) \\ u_2(t) \\ u_3(t)
    \end{bmatrix} =
    \begin{bmatrix}
      (3+\cos(6\pi t/50))\cos(4\pi t/50) \\ (3+ \cos(6\pi t/50))\sin(4\pi t/50) \\ \sin(6\pi t/50) 
    \end{bmatrix},
  \end{equation*}
  which we extend in a fixed direction to obtain nodal function values that we use a preliminary initial data.
  For triangulations $\Th = \cT_k$ of $\o$ into triangles given by halved squares with side length $\hat{h} = 2^{-k}$, we thus define the preliminary initial data via
  \begin{equation*}
      \ty_h^0(z) = \begin{bmatrix} u_1(z_1) \\ u_2(z_1) \\ u_3(z_1)+z_2 \end{bmatrix}
    , \quad
    \nabla \ty_h^0(z) = \begin{bmatrix} \eta(z_1) \dot{u}_1(z_1) & 0 \\ \eta(z_1) \dot{u}_2(z_1) & 0 \\ 0 & 1 \end{bmatrix},
  \end{equation*}
  with $\eta(t) = (\dot{u}_1(t)^2 + \dot{u}_2(t)^2)^{-2}$ for every vertex $z=(z_1,z_2) \in \Nh$. As before, the actual initial data $y_h^0$ is obtained from a relaxation of $\ty_h^0$, cf. Example~\ref{ex:twist180}. 
  The iteration numbers, final energies, tangent-point potentials and isometry errors obtained with exponent $q = 5$ on a sequence of triangulations, $k=1,2,3$, with step size $\tau \sim \hat{h}$ and $\TP$ parameter $\rho \sim \hat{h}$ are shown in Table~\ref{tab:trefoil-table}. 
  The numerical equilibrium configuration for the finest triangulation $\cT_3$ is depicted in Figure~\ref{fig:trefoil}, which also contains an illustration of the pseudo force $f_\tp$ induced by the potential $\TP$ at the final iterate. 
  It is clearly visible that the pseudo force acts in an almost normal direction whereas the ``natural direction'' of a self-avoiding pseudo force in this example should be tangential to the plate and normal to its boundary at almost-contact points.
  The experimental order of convergence in $L^2$ for the piecewise linear interpolants of the discrete solutions in this example is approximately~$\mathcal{O}(h^{0.44})$.
  Note that the expected contact zone in this example is entirely a subset of the domain boundary $\p \o$. 
  The singular nature (``edge-to-edge'') of the almost-contact requires a small time step size to avoid energy blowups and a careful choice of the $\TP$ parameter $\rho$, which has do be chosen large enough to guarantee self-avoidance, but at the same time should be chosen as small as possible to minimize potential repulsive effects other than the prevention of self-intersections. 
  In this situation, additionally including a tangent-point potential of the (one-dimensional) boundary curve of the domain in the energy functional -- or even replacing the surface potential with the boundary curve potential -- might provide a remedy.
\end{example}

\subsection{O-shaped bilayer plate}
\begin{figure}
  \centering
  \includegraphics[width=0.98\textwidth]{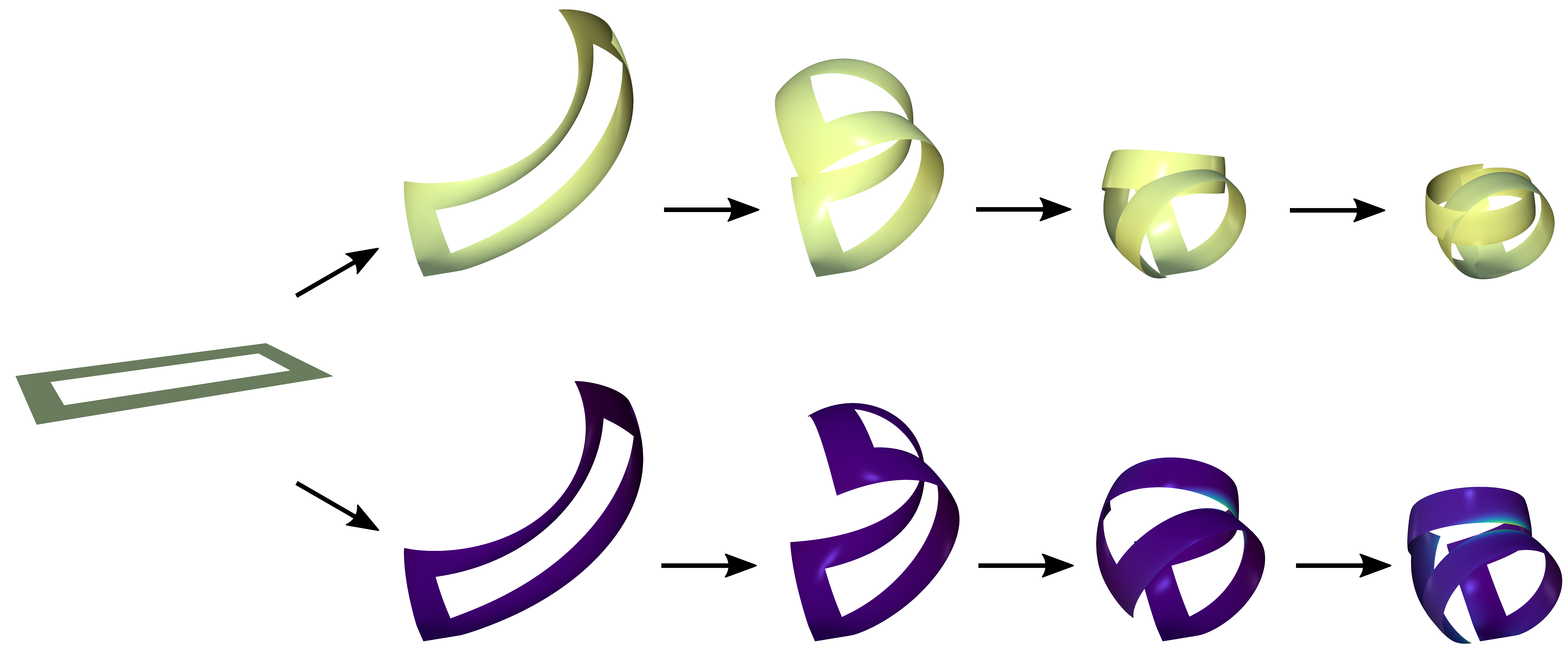}
  \put(-345,119){$\rho = 0$}
  \put(-345,12){$\rho > 0$}
  \caption{Snapshots of the evolution in Example~\ref{ex:oshape} for the triangulation $\cT_4$ after 0, 70, 450, 1050 and 3077 iterations (bottom) in comparison with a similar evolution corresponding to the bilayer energy without $\TP$ term (top). \label{fig:oshape_evo}}
\end{figure}

In order to investigate the effect of including the tangent-point potential in the bilayer plate model we consider an O-shaped bilayer plate which is horizontally clamped on one of its corners. We choose a material mismatch that guarantees self-intersections if $\rho=0$, i.e. if the tangent-point potential is neglected, cf. Figure~\ref{fig:oshape_evo}.
\begin{table}
\centering
\begin{tabular}{c c c c c c c}\hline
  $k$ & $\tau = \rho$ & $N_\mathrm{iter}$ & $E_h[y_h^\infty]$ & $\TP_{h}[y_h^\infty]$ & $\delta_\mathrm{iso}[y_h^\infty]$ \\ \hline
  $2$ & \num{0.2500}   & 1229 & \num{-7.59613}  & \num{0.567855}  & \num{2.92514} \\
  $3$ & \num{0.1250}   & 1563 & \num{-0.182432} & \num{0.312536}  & \num{1.08479} \\
  $4$ & \num{0.0625}   & 3077 & \num{ 2.94449}  & \num{0.231738}  & \num{0.439736} \\
\end{tabular}
  \caption{Iteration numbers, total energy, tangent-point potential and isometry error of the final iterates in the O-shaped plate Example~\ref{ex:oshape}.}
  \label{tab:oshape_table}
\end{table}

\begin{example}\label{ex:oshape}
  Let $\alpha= 0.75$ in the bilayer energy and $f = [0,0,0]^\top$ in 
  \[
  \o = (-5,5) \times (-2,2) \setminus [-4,4] \times [-1,1]
  \]
  with clamped boundary conditions along the corner 
  \[
  \GD = \{-5\} \times [-2,-1] \cup [-5,-4] \times \{-2\},
  \]
  i.\,e. we have $y_\DD = [x,0]^\top$ and $\phi_\DD = [I_2,0]^\top$.
  We compute the numerical equilibrium states for a sequence of triangulations $\Th = \cT_k$, $k=2,3,4$, consisting of halved squares with side length $\hat{h} = 2^{-k}$ amounting to 12288 triangles and 58455 degrees of freedom in the case $k=4$. 
  We choose the step size and TP parameter as $\tau = \rho = \hat{h}$ and use the exponent $q = 5$. The resulting iteration numbers, final energies, tangent-point potentials and isometry errors are shown in Table~\ref{tab:oshape_table}. 
  In this example, the experimental order of convergence in $L^2$ for the piecewise linear interpolants of the discrete solutions is approximately~$\mathcal{O}(h^{1.3})$.
  The evolution corresponding to the discrete gradient flow for the triangulation $\cT_4$ is shown in Figure~\ref{fig:oshape_evo}, in comparison with an evolution corresponding to a discrete gradient flow for the case~$\rho=0$. 
  In all considered cases self-intersections were prevented until the numerical equilibrium configuration was reached.
\end{example}

\subsection{Self-coiling bilayer plate}
To conclude our examples we consider two rectangular bilayer plates with different lengths which are clamped horizontally on one short side. Each plate can be expected to roll itself up into a cylindrical shape as a consequence of the material mismatch. Indeed, for the considered boundary conditions, analytical minimizers of the bilayer energy $E_\mathrm{bil}$ are given by parametrizations of cylinders with radius $\a^{-1}$, see~\cite{Schmidt07}. Hence, analytical minimizers cannot be injections if the length of the long side of the rectangular domain is greater than $2\pi\a^{-1}$.
\begin{figure}
  \centering
  \subfloat{\includegraphics[height=3.4cm]{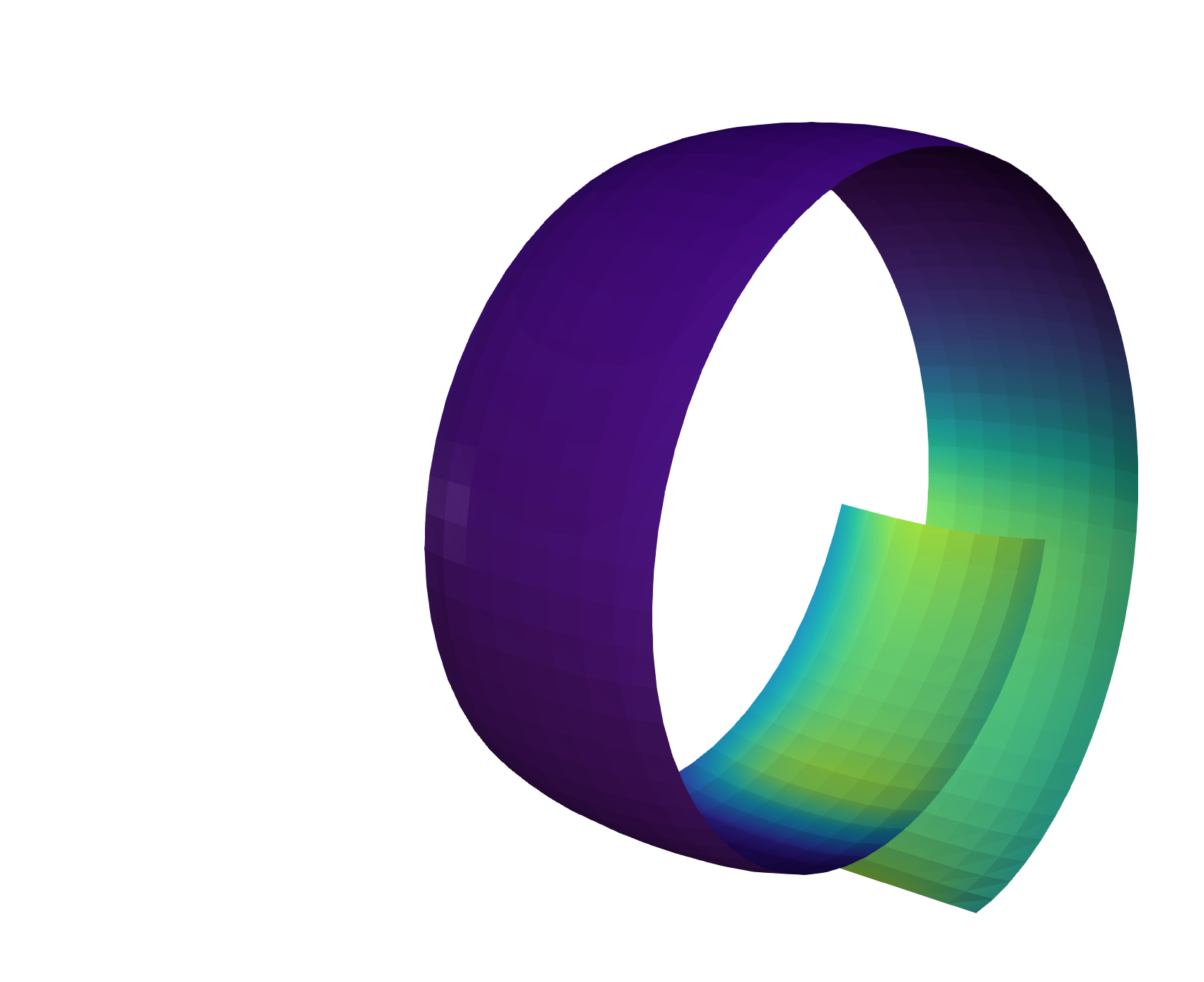} \put(-80,0){$\o_\mathrm{short}$} \hspace{1.2cm}}
  \subfloat{\includegraphics[height=3.4cm]{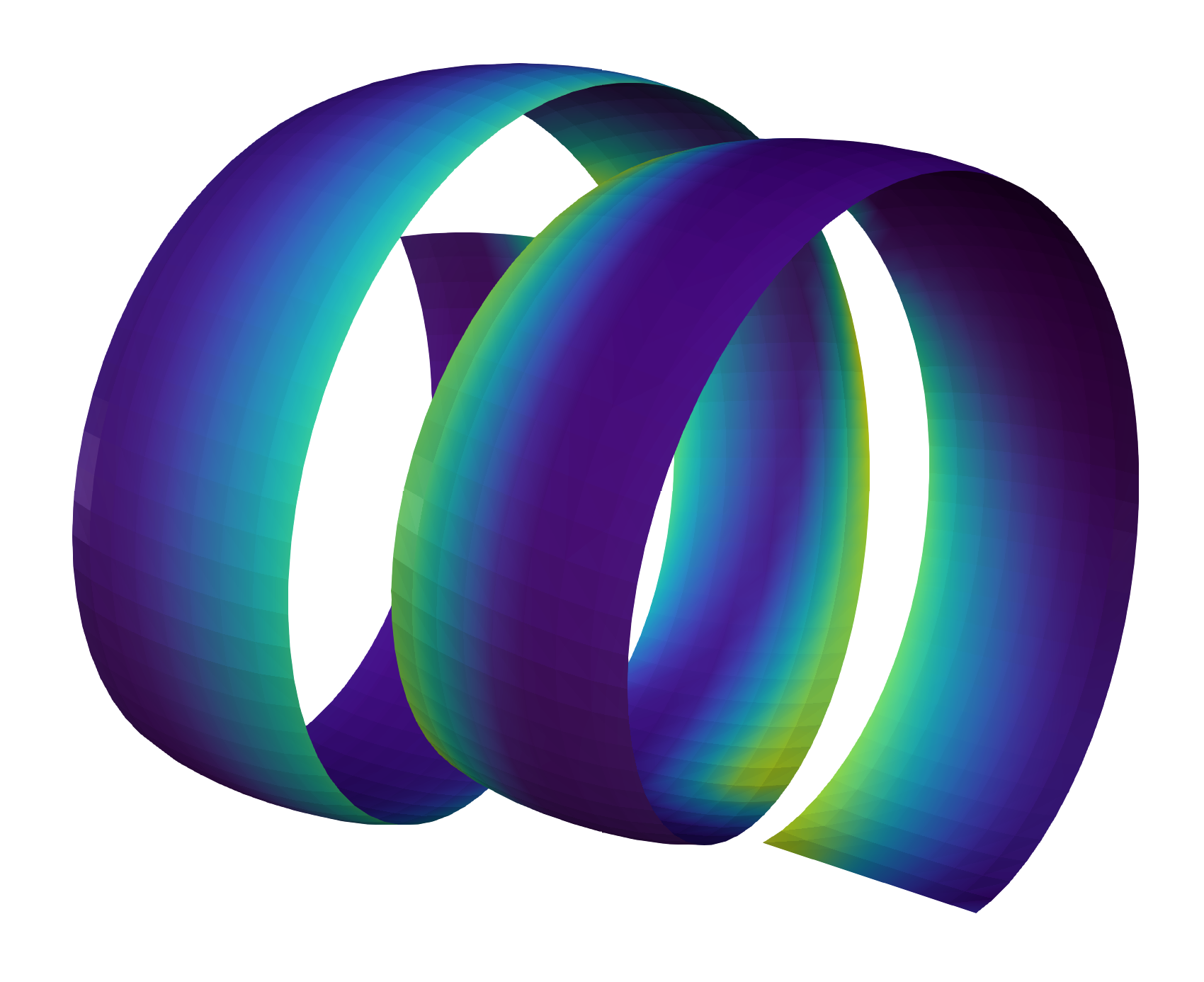} \put(-80,0){$\o_\mathrm{long}$} \hspace{1.5cm}}
  \caption{Numerical equilibrium configurations of the self-coiling bilayer plates in Example~\ref{ex:bilayer}. Self-repulsive forces in the longer plate become leading to a stationary configuration that is reminiscent of a corkscrew. \label{fig:bilayer}}
\end{figure}
\begin{table}
  \begin{tabular}{c c c c c c c c}\hline
  Domain & $k$ & $\tau$ & $\rho$ & $N_\mathrm{iter}$ & $E_h[y_h^\infty]$ & $\TP_{h}[y_h^\infty]$ & $\delta_\mathrm{iso}[y_h^\infty]$ \\ \hline
  $\o_\mathrm{short}$ & $2$ & \num{0.0125} & \num{0.000625} &  22118 & \num{1.03537} & \num{14.2236}  & \num{0.14914} \\
  $\o_\mathrm{long}$  & $2$ & \num{0.0125} & \num{0.000625} & 104683 & \num{1.30437} & \num{21.3715}  & \num{0.325157} \\
  $\o_\mathrm{short}$ & $3$ & \num{0.00625} & \num{0.0003125} &  46129 & \num{2.09189} & \num{13.3693}  & \num{0.0710055} \\
  $\o_\mathrm{long}$  & $3$ & \num{0.00625} & \num{0.0003125} & 199162 & \num{3.87973} & \num{21.7176}  & \num{0.155099} \\
\end{tabular}
\caption{Iteration numbers, total energies, tangent-point potentials and isometry errors in the self-coiling bilayer Example~\ref{ex:bilayer} for different refinement levels~$k$.}
\label{tab:bilayer_table}
\end{table}

\begin{example}\label{ex:bilayer}
  We consider the minimization problem for the bilayer energy with material mismatch $\a = 1$ for two plates with reference configurations 
  \[
  \o_\mathrm{short} = (0,10) \times (0,1), \quad \o_\mathrm{long} = (0,20) \times (0,1).
  \] 
  We let $f = [0,0,0]^\top$ and prescribe clamped boundary conditions $y_\DD = [x,0]^\top$ and $\phi_\DD = [I_2,0]^\top$ along the edge $\GD = (0) \times (0,1)$. We compute the discrete gradient flow on triangulations $\Th = \cT_k$, $k=2,3$, consisting of halved squares with side length $\hat{h} = 2^{-k}$, using step size $\tau = \hat{h}/20$ and TP parameter $\rho = \hat{h}/400$.
  The TP exponent is chosen as $q=8$ in this example. For the refinement level $k=3$ the discrete gradient flows terminate after $46129$ and $199162$ iterations for $\o_\mathrm{short}$ and $\o_\mathrm{long}$, respectively, at configurations shown in Figure~\ref{fig:bilayer}. Corresponding final energies, values of TP and isometry errors are listed in Table~\ref{tab:bilayer_table}. The stationary configurations demonstrate the effect of the tangent-point potential as a self-repulsive force. The material mismatch in the short plate leads to an almost one-dimensional deformation into a spiral shape with an incircle radius of about $1.2$ which is comparable to the radius $1$ of the analytical minimizer for the pure bilayer energy in this setting. Self-intersections are successfully prevented but the \emph{distance from self-contact} is noticeable for the chosen parameters. The long plate, on the other hand, undergoes an additional out-of-plane deformation along the $y$-axis and develops a corkscrew-like configuration with a similar approximate incircle radius. Self-intersections are successfully prevented, but the effects of the self-repulsive force are noticeable. 
\end{example}

\section{Conclusion}

The discretization of the two-dimensional tangent-point potential can be employed to avoid self-intersections in the simulation of bending isometries. The proposed semi-implicit discrete gradient flow method leads to linear problems in every time step and is practical in this regard. The computations involved in the assembly of the linear problems are expensive and special care should be devoted to ensure an efficient implementation. Depending on the specific problem under consideration a careful choice of the involved parameters as well as the choice of a small (pseudo-)time step size might be necessary and can lead to many iterations, limiting the practical efficiency of the method. This is particularly evident for problems in which configurations in a neighborhood of a singular point of the tangent-point potential might occur, e.\,g. almost-edge-to-edge contact with almost parallel tangent planes. A possible remedy for such difficulties is the use of an augmented tangent-point potential composed of the two-dimensional surface potential as well as the one-dimensional tangent-point potential of the boundary curve of the plate.

For the class of problems with non-singular almost-contact region our findings indicate that the proposed method is a practical way of avoiding self-intersections. The two-dimensional self-avoiding plate model resulting from the inclusion of the tangent-point potential does, however, not consider physical contact phenomena such as friction, and a physical justification or rigorous derivation remain open.

\smallskip
\noindent {\small{\textbf{Acknowledgements} The authors gratefully acknowledge the support by the Deutsche Forschungs\-gemeinschaft in the Research Unit~3013 \emph{Vector- and Tensor-Valued Surface PDEs} within the sub-project \emph{TP4: Bending plates of nematic liquid crystal elastomers}.}}

\bibliographystyle{siam}
\bibliography{references}

\end{document}